\documentclass[12pt,a4paper,reqno]{amsart}

\usepackage{hyperref}
\usepackage{mathrsfs}
\usepackage{amstext}
\usepackage{amscd}
\usepackage{latexsym}
\usepackage{amsfonts}
\usepackage{amssymb}
\usepackage{amsmath}
\usepackage{amsthm}
\usepackage{color}
\usepackage{multicol}
\usepackage[pdftex]{graphicx}
\usepackage{url}

\usepackage[russian,english]{babel}
\usepackage[cp1251]{inputenc}

\theoremstyle{plain}
\newtheorem{theorem}{Theorem}[section]

\newtheorem{lemma}[theorem]{Lemma}
\newtheorem{propos}[theorem]{Proposition}
\newtheorem{corol}[theorem]{Corollary}

\theoremstyle{definition}
\newtheorem{definition}[theorem]{Definition}

\newtheorem{conjecture}[theorem]{Conjecture}

\newtheorem{remark}[theorem]{Remark}
\newtheorem{example}[theorem]{Example}
\newtheorem{notation}[theorem]{Notation}

\newtheorem{letteredtheorem}{Theorem}

\theoremstyle{thmstyletwo}

\begin{document}

\author[O.~Katkova]{Olga Katkova}

\address{Department of Mathematics, University of Massachusetts Boston}
\email{Olga.Katkova@umb.edu}

\author[M.~Tyaglov]{Mikhail Tyaglov}

\address{Moscow Center for Fundamental and Applied Mathematics, 119991, Moscow, Russia} 
\address{School of Mathematical Sciences and MOE-LSC, Shanghai Jiao Tong University}
\email{tyaglov@gmail.com,\ \ tyaglov@sjtu.edu.cn}

\author[A.~Vishnyakova]{Anna Vishnyakova}

\address{Department of Mathematics,  Holon Institute of Technology, 
Holon, Israel}
\email{annalyticity@gmail.com}

\title[A generalization of the Hawaii conjecture]{On the number of real zeroes 
of a homogeneous differential polynomial and a generalization of the Hawaii 
conjecture}

\date{\small \today}



\begin{abstract}
For a given real polynomial $p$ we study the possible number of real roots of 
a differential polynomial $H_{\varkappa}[p](x) =
\varkappa\left(p'(x)\right)^2-p(x)p''(x),  \varkappa \in \mathbb{R}.$  In the special case 
when all real zeros of the polynomial $p$ are simple, and all roots of its derivative $p'$ 
are real and simple, the distribution of zeros of $H_{\varkappa}[p]$ is completely described 
for each real $\varkappa.$ We also provide counterexamples to 
two Boris Shapiro's conjectures about the number of zeros of the function 
$H_{\frac{n-1}{n}}[p].$

\end{abstract}

\maketitle

\setcounter{equation}{0}

\section{Introduction}\label{section:intro}

\setcounter{equation}{0}

\medskip
\begin{center}
\textit{To our dear teacher  and friend Victor Katsnelson,
an outstanding mathematician and a bright person, with love and gratitude. 
}
\end{center}
\medskip

In this paper, we investigate the possible number of real zeros of
special homogeneous differential polynomials.

\begin{notation}\label{Not:Q.function}
Let $p$ be a real polynomial and $\varkappa$ be a complex number.
We  denote by $Q_{\varkappa}[p]$ the following function \it{associated with} $p$
\begin{equation}\label{main.function.2}
Q_{\varkappa}[p](x) = p^{1-\varkappa}(x)\cdot\dfrac{d}{dx}
\left(\dfrac{p^{\varkappa}(x)}{p'(x)}\right)=
\dfrac{\varkappa\left(p'(x)\right)^2-p(x)p''(x)}{\left(p'(x)\right)^2}.
\end{equation}
By  $H_{\varkappa}[p]$ we  denote the following function:
\begin{equation}\label{main.function.2.numerator}
H_{\varkappa}[p](x)\stackrel{def}{=}
\varkappa\left(p'(x)\right)^2-p(x)p''(x)=\left(Q_{\varkappa}[p](x)\right)(p'(x))^2.
\end{equation}
\end{notation}

We study the number of real roots of $Q_{\varkappa}[p]$ and $H_{\varkappa}[p]$,
and we need the following notation.

\begin{notation}\label{notation.number.of.zeroes}
For a real polynomial $p$, by $Z_\mathbb{C}(p)$ we denote
the number of non-real zeroes of $p$ counting multiplicities. If $f$ is a real 
rational function, then $Z_\mathbb{R}(f)$ will denote the number of real zeroes 
of $f$ counting multiplicities. Generally, the number of zeroes of a real 
rational function $f$ on a set $X \subset\mathbb{R}$ counting multiplicities 
will be denoted by $Z_X(f)$. In particular, $Z_\mathbb{R}(f)=
Z_{(-\infty,+\infty)}(f).$
\end{notation}

Our interest is concentrated on bounding the number of real
zeroes of $Q_{\varkappa}[p]$ and $H_{\varkappa}[p]$, that is on 
bounding $Z_\mathbb{R} (Q_{\varkappa}[p])$ and 
$Z_\mathbb{R} (H_{\varkappa}[p])$ for $\varkappa\in\mathbb{R}.$

Note that if a polynomial $p$ doesn't have multiple zeros, then
\begin{equation}\label{main.functions.zeros}
Z_\mathbb{R} (Q_{\varkappa}[p])=Z_\mathbb{R} \left(H_{\varkappa}[p]\right).
\end{equation}

The following statement belongs to E.~Laguerre.
\begin{letteredtheorem} \label{Theorem:Laguerre Inequality}
If a real polynomial $p$ has only  real and simple zeros,  then
$H_1[p](x)=\left(p'(x)\right)^2-p(x)p''(x)>0$ for all $x\in \mathbb{R}.$
\end{letteredtheorem}

The Laguerre inequality plays an important role in the study of distribution of zeros of 
real entire functions and in the understanding the nature of the Riemann $\xi-$function 
(see~\cite{CravenCsordas1}, \cite{CravenCsordas2}, 
\cite{CsordasEscassut}--\cite{CsordasVishnyakova} and the references therein for the 
generalization of the Laguerre inequality, as well).

The following generalization of the Laguerre inequality is known as the Hawaii conjecture.

\begin{conjecture}[the Hawaii Conjecture, \cite{CCS_1987}]\label{Theorem:Hawaii Conjecture}
 If all real roots of a real polynomial $p$ are simple,  then 
\begin{equation}\label{Hawaii}
Z_\mathbb{R} \left(H_1[p]\right)\leq Z_\mathbb{C} (p).
\end{equation}
\end{conjecture}

The Hawaii conjecture was posed by T.~Craven, G.~Csordas, and W.~Smith 
in~\cite{CCS_1987} in 1987, and was proved by Mikhail Tyaglov in 2011 
(\cite{Tyaglov_Hawaii}).

A different situation when the polynomials $H_{\varkappa}[p]$ appear was discovered by 
J.L.W.V.~Jensen in~\cite{Jensen}. Consider the function
\begin{equation}\label{Jensen_function}
\Phi_p (x,y) =\left|p(x+iy)\right|^2.
\end{equation}
The function $\Phi_ p(x,y)$ is real-analytic non-negative function in $(x,y)$, whose zeros are 
zeros of $p(z)=p(x+iy).$ Assume that $\deg p=n$ and expand $\Phi_ p(x,y)$ in the variable 
$y$. We obtain
\begin{equation}\label{Jensen_expansion}
\Phi_p (x,y) =\sum_{k=0}^n P_k(x) \frac{y^{2k}}{(2k)!},
\end{equation}
where
\begin{equation}\label{polynomialPk}
P_k(x) =\sum_{j=0}^{k }(-1)^{j+k}\binom{2k}{j}p^{(j)}(x)p^{(2k-j)}(x).
\end{equation}
In particular, $P_0(x)=p^2(x)$ and $P_1(x)=\left(p' (x)\right)^2-p(x) p''(x).$
Jensen obtained the following sufficient and necessary conditions for a polynomial to have 
only real  simple roots: {\it the polynomials $P_k (x)$ are strictly positive for all 
$k= 1, 2, \ldots, n.$} While studying unpublished notes left after Jensen's death in 1925, 
G.~P\'{o}lya discovered a different criterion of real-rootedness.

\begin{letteredtheorem}[G.~P\'{o}lya]\label{Polyarealroot}
 A real polynomial $p$ of degree $n$ has only  real  simple zeros if and only if the  polynomials
$$G_k[p](x)=(n-k)\left(p^{(k)}(x)\right)^2-(n-k+1)p^{(k-1)}(x)p^{(k+1)}(x)$$ 
are strictly positive for all $k=1,2,\ldots, n.$
\end{letteredtheorem}

Inspired by the Hawaii conjecture and motivated by these two criteria, Boris~Shapiro formulated 
the following two conjectures.

\begin{conjecture}[Boris~Shapiro, \cite{Shapiro}]\label{Shapiro1}
 For any real polynomial $p$ of degree $n\geq 2$ all whose real roots are simple
\begin{equation}\label{Sh1}
Z_\mathbb{R} \left[(n-1)\left(p'\right)^2-n p  p''\right]\leq Z_\mathbb{C} \left(p\right).
\end{equation}
If this fact was true, the Hawaii conjecture would hold for $G_1.$
\end{conjecture}

\begin{conjecture}[Boris~Shapiro, \cite{Shapiro}]\label{Shapiro2}
For any real polynomial $p$ of even degree $n\geq 2$,
\begin{equation}\label{Sh2}
Z_\mathbb{R} \left[(n-1)\left(p'\right)^2-n p  p''\right]+ Z_\mathbb{R}(p)>0.
\end{equation}
The latter inequality is trivially satisfied for odd degree polynomials.
\end{conjecture}

\setcounter{equation}{0}

\section{Main results}\label{section:lower.bound}

It was observed in~\cite{Tyaglov.Atia.2021} that the value $\varkappa=\dfrac{n-1}{n}$ is of special 
significance for functions $Q_{\varkappa}[p].$ Although for polynomials all whose zeroes are real  
the properties of $Q_{\tfrac{n-1}{n}}[p]$ are close to the ones of $Q_{1}[p],$  there are polynomials 
of degree $n\geqslant4$ with non-real zeros not satisfying  Conjecture~\ref{Shapiro1}. 
In~\cite{Tyaglov.Atia.2021} there was given the example of such a polynomial of degree $4$ 
having two real and two non-real zeros.
\begin{letteredtheorem}[M. Tyaglov, M.J. Atia, \cite{Tyaglov.Atia.2021}]\label{counterTA}
 If
\begin{equation}\label{counter1}
p(x)=(x^2+a^2)(x+a^2)(x-1),\qquad a\in\mathbb{R}\backslash\{-1,0,1\},
\end{equation}
 then
$$Z_\mathbb{R} \left(H_{\frac{3}{4}}[p]\right)=4.$$
\end{letteredtheorem}

It is obvious that $Z_\mathbb{C}(p)=2<4.$ Therefore, Theorem~\ref{counterTA} gives a 
counterexample to the Conjecture~\ref{Shapiro1}. 

However, the question still remains: whether or not counterexamples given by polynomials of higher degrees exist.
In the present paper, for each $n=5, 6,\ldots,$ we construct a real polynomial $p$ of degree $n$ that doesn't 
satisfy the Conjecture~\ref{Shapiro1}.  
\begin{propos}\label{counterAnna}
 If $n\geq 5$ and
\begin{equation}\label{counter2}
p(x)=(x-1)^n+(x+1)^n,
\end{equation}
 then
$$Z_\mathbb{R} \left(H_{\frac{n-1}{n}}[p]\right)=2n-4,$$
 while
$$Z_\mathbb{C}(p)=2\left\lfloor\dfrac n2\right\rfloor < 
Z_\mathbb{R} \left(H_{\frac{n-1}{n}}[p]\right).      $$
\end{propos}
It follows from Theorem~\ref{counterTA} and Proposition~\ref{counterAnna} that for 
$\varkappa=\frac{n-1}{n}, \ \ n=4,5,6,\ldots,$ the generalization of the Hawaii conjecture 
\begin{equation}\label{Hawaiigeneral}
Z_\mathbb{C} \left(Q_{\varkappa}[p]\right)\leq Z_\mathbb{R} (p)
\end{equation}
is not valid.

There arises a natural question: for which values of $\varkappa$ the Hawaii conjecture 
continues to be true?

We investigate the number of real zeros of the function $Q_{\varkappa}[p], \ \varkappa \in \mathbb{R},$
in the case when the derivative $p'$ has {\it only real simple zeros.} The statement bellow can be considered 
as a generalization of the Hawaii conjecture in this case. 
\begin{theorem}\label{Hawaiitrue}
Let $\varkappa\geq\frac{n-1}{n}.$  For any real polynomial $p $ of degree $n\geq 2$ 
all whose real zeros are simple, and all zeros of whose derivative $p'$ are real and simple,
\begin{equation}
\label{t1}
Z_{\mathbb{R}}\left(H_{\varkappa}[p]\right)=Z_{\mathbb{C}} (p). 
\end{equation}
\end{theorem}

The following two statements give estimations for the number 
$Z_{\mathbb{R}}\left(H_{\varkappa}[p]\right)$ 
when $\varkappa<\frac{1}{2}.$
\begin{theorem}\label{Theorem2}
 Let $\varkappa\leq 0.$  For any real polynomial $p $ of degree $n\geq 2$ all whose real 
 zeros are simple, and all zeros of whose derivative $p'$ are real and simple,
\begin{equation}
\label{t2}
Z_{\mathbb{R}} \left(H_{\varkappa}[p]\right)=n+Z_{\mathbb{R}} (p) -2.
\end{equation}
\end{theorem}

\begin{theorem}\label{Theorem3}
Let $0<\varkappa <1/2.$  For any real polynomial $p $ of degree $n\geq 2$ all whose real zeros are 
simple, and all zeros of whose derivative $p'$ are real and simple,
\begin{equation}
\label{t3}
Z_{\mathbb{C}}(p)-2 \leq Z_{\mathbb{R}} \left(H_{\varkappa}[p]\right)\leq n+Z_{\mathbb{R}} (p) -2.
\end{equation}
\end{theorem}

We will show that both estimations in Theorem \ref{Theorem3} are sharp. In particular, we prove 
the following theorem.
\begin{theorem}\label{Theorem3sharp}
For each $\varkappa,   0<\varkappa<\frac{1}{2},$ and each $n\geq 3$ there exists a polynomial $p_n$ 
of degree $n$ all whose real zeros are simple, and all zeros of whose derivative $p'$ are real 
and simple, such that 
$$Z_{\mathbb{R}} \left(H_{\varkappa}[p_n]\right)= n+Z_{\mathbb{R}} (p_n) -2.$$
\end{theorem}
Since 
$$ n+Z_{\mathbb{R}} (p) -2\leq Z_{\mathbb{C}} (p) \  \Leftrightarrow 
\ Z_{\mathbb{R}} (p) \leq 1,$$  Theorem \ref{Theorem2} and Theorem \ref{Theorem3sharp} 
say that, in general, the Hawaii conjecture
is not valid for $\varkappa<\frac{1}{2}.$

In~\cite{Tyaglov.Atia.2021} M.~Tyaglov and M.J.~Atia found the upper and lower bounds 
of $Z_{\mathbb{R}} \left(H_{\varkappa}[p]\right) $ for all real $\varkappa$ when a polynomial 
$p$ has only real zeros. Their result can be reformulated as follows.
\begin{letteredtheorem}[M.~Tyaglov, M.J.~Atia, \cite{Tyaglov.Atia.2021}]\label{TheoremTA}
Let $\frac{k-1}{k}\leq \varkappa <\frac{k}{k+1}, \ k=2,3,\ldots,n-1.$  If a real polynomial $p$ of degree 
$n\geq 2$  has only real and simple zeros, then
\begin{equation}
\label{tTA}
2 \leq Z_{\mathbb{R}} \left(H_{\varkappa}[p]\right)\leq 2n-2k.
\end{equation}
\end{letteredtheorem}

We will show that all estimations in Theorem~\ref{TheoremTA} are sharp.

Analyzing the statements of Theorems \ref{Hawaiitrue} -- \ref{Theorem3sharp} and 
Theorem\ref{TheoremTA}, one can guess that for any real polynomial $p$ of degree 
$n$ the Hawaii conjecture (\ref{Hawaiigeneral}) holds when $\varkappa>\frac{n-1}{n}.$ 
The statement below, which is easy to check  by direct calculations,  shows that this 
conjecture is not true.

\begin{propos}\label{counterOlga}
Let $n\geq 3,$ and 
\begin{equation}\label{counterO}
p(x)=x^n+ax^{n-2}, \ \, a>0.
\end{equation}
Then for all $\varkappa \in \left[\frac{n-1}{n}, \ \frac{(2n-3)^2}{4n(n-2)}\right]$ 
\begin{equation}\label{counterOO}
Z_\mathbb{R} \left(Q_{\varkappa}[p]\right)=4.
\end{equation}
\end{propos}

\begin{remark}\label{counterOlga1}
Note that $\frac{(2n-3)^2}{4n(n-2)}>\frac{n-1}{n}.$  So, in general, the number $\frac{n-1}{n}$ 
cannot serve as the infimum for the set of all real $\varkappa$ for which (\ref{Hawaiigeneral}) is valid.
\end{remark}

\begin{remark}\label{counterOlga2}
It is worth noting that the polynomial in the above example has $0$ as a root of multiplicity $(n-2).$ 
Any polynomial of the form $p_{2n, a, C}(x)=x^{2n}+ax^{2n-2}+C, \ \ C\neq 0,$ has either two or zero 
real roots, and any polynomial $p_{2n+1, a, C}(x)=x^{2n+1}+ax^{2n-1}+C, \ \ C\neq 0,$ has only one real root. 
So, if $\deg p \geq 5,$ even a small shift of the polynomial in Proposition~\ref{counterOlga} gives 
$Z_\mathbb{C}(p)\geq 4,$ which makes (\ref{counterOO}) consistent with  (\ref{Hawaiigeneral}). 
Probably, the question about the infimum for the set of all real $\varkappa$ for which (\ref{Hawaiigeneral}) 
is valid will make more sense for polynomials, all whose real roots are simple. 
\end{remark}

Unfortunately, in the general case, if there aren't any restrictions on roots of the derivative $p'$ 
of a real polynomial $p,$ we could find only a lower bound for the number of real zeroes of the 
function $Q_{\varkappa}[p]$ for an arbitrary real 
polynomial $p$ whenever $\varkappa\in\mathbb{R}$. 

\begin{notation}\label{notation.number.of.zeroes.new}
We denote by $\mathring{Z}_{X}(f)$ the number of zeros of the function 
$f$ on the set $X$ without counting their multiplicities.
\end{notation}

\begin{theorem}\label{Theorem:lower.bound}
Let $p$ be a real polynomial of degree $n\geqslant2$.  The following inequalities hold
\begin{equation}\label{main.result.lower.bound.1}
Z_\mathbb{R}\left(Q_{\varkappa}[p]\right)\geqslant \mathring{Z}_\mathbb{R}(p')+1-
\mathring{ Z}_\mathbb{R}(p),\qquad\text{for}\qquad\varkappa > \dfrac{n-1}{n},
\end{equation}
\begin{equation}\label{main.result.lower.bound.2}
Z_\mathbb{R}\left(Q_{\varkappa}[p]\right)\geqslant \mathring{Z}_\mathbb{R}(p')-1- 
\mathring{Z}_\mathbb{R}(p),\qquad\text{for}\qquad 0<\varkappa\leqslant\dfrac{n-1}{n}.
\end{equation}
If $p(x)\neq (x-\alpha)^n,$ then
\begin{equation}\label{main.result.lower.bound.3}
Z_\mathbb{R}\left(Q_{\varkappa}[p]\right)\geqslant \mathring{Z}_\mathbb{R}(p')-1+ 
\mathring{Z}_\mathbb{R}(p),\qquad\text{for}\qquad \varkappa\leqslant 0.
\end{equation}
\end{theorem}

Now let us return to Boris~Shapiro's Conjecture~\ref{Shapiro2}. In their recent 
preprint \cite{MaMa}, Lande Ma and Zhaokun Ma discussed this conjecture for entire 
transcendental functions.  For any even number $n\geq 4$ we  give an example of a 
real polynomial $p$ of degree $n$ that doesn't satisfy Conjecture~\ref{Shapiro2} .

\begin{propos}\label{counterOlga2}
 Let $n\geq 2,$ and 
\begin{equation}\label{counterO2}
p_{2n}(x)=\frac{x^{2n}}{2n}+\frac{x^2}{2}+1.
\end{equation}
Then 
$$Z_\mathbb{R} \left[(2n-1)(p'(x))^2-2n p(x)p''(x) \right]=0.$$
\end{propos}

The following theorem describes a sufficiently broad class of polynomials satisfying 
Conjecture~\ref{Shapiro2}.

\begin{theorem}\label{conjectureSh2}
 Let $p $ be a real polynomial of even degree $n\geq 4.$ Assume that there 
 exists $k=1,2,\ldots, n-2,$ such that all the zeros of the derivative $p^{(k)}$ 
are real. Then
\begin{equation}
\label{Sh2}
Z_{\mathbb{R}}\left(H_{\frac{n-1}{n}}[p]\right)+Z_{\mathbb{R}} (p) >0.
\end{equation}
\end{theorem}

The paper has the following structure. In Section~\ref{Sec3} we discuss counterexamples to 
Boris~Shapiro's  conjectures. In Section~\ref{Sec4} we introduce a function $M[p],$ which 
plays a major role in all our constructions. As an example of applications of $M[p],$ in this 
Section we prove two theorems: Theorem~\ref{Theorem:lower.bound} on a lower bound for 
the number $Z_\mathbb{R} \left(Q_{\varkappa}[p]\right),$ and Theorem~\ref{conjectureSh2} 
describing polynomials that satisfy Boris Shapiro's conjecture~\ref{Shapiro2}. In the rest of our 
paper we investigate the zero distribution of the functions $Q_{\varkappa}[p]$, where $p$ are 
real polynomials, all whose real zeros are simple and whose derivative $p'$ has only real and simple 
roots. Sections~\ref{Sec5}--\ref{Sec9}  are devoted to the case $\varkappa>\frac{n-1}{n}.$  We 
develop an idea we found in \cite[p.106, Problem 725]{FaddeevSominskii}  and  prove that in this 
case the Hawaii conjecture is true. The result from \cite{FaddeevSominskii} relates to the case 
when $\varkappa=0$ and all roots of $p$ are real and simple. In Section~\ref{Sec10} we obtain 
a result (Lemma~\ref{lemma8}) similar to Rolle's theorem that gives a relation between roots of 
$H_{\varkappa}[p]$ and $H_{2-\frac{1}{\varkappa}}[p'].$ This Lemma allows to investigate the 
zero distribution of the function $Q_{\varkappa}[p]$ in the case $-\infty<\varkappa<\frac{n-1}{n}$ 
(Sections~\ref{Sec11}--\ref{Sec13}). In Section~\ref{Sec14} we discuss the accuracy of estimates in 
Theorem~\ref{Theorem3}. Sections~\ref{Sec15} contains some additional results about the zero 
distribution of $Q_{\varkappa}[p]$ in the case when $\frac{1}{2} \leq \varkappa<\frac{n-1}{n}$ and 
all roots of $p$ are real. In Section~\ref{Sec16} we describe possible distribution of roots of 
$Q_{\varkappa}[p]$ when $\frac{1}{2} \leq \varkappa<\frac{n-1}{n}$ in the general case.

\section{ Counterexamples to Boris~Shapiro's Conjectures: Proof of Theorem~\ref{counterTA}, 
and Propositions~\ref{counterAnna} and~\ref{counterOlga2} }  \label{Sec3}

The fact that Conjecture~\ref{Shapiro1} is not true for $n =\deg p=4$ was established 
in~\cite{Tyaglov.Atia.2021}. We provide a proof of Theorem~\ref{counterTA} 
here for completeness and convenience of readers.

{\bf Proof of Theorem~\ref{counterTA}.} The polynomial of the fourth degree
\begin{equation}\label{counter1}
p(x)=(x^2+a^2)(x+a^2)(x-1),\qquad a\in\mathbb{R}\backslash\{-1,0,1\},
\end{equation}
has two distinct real zeroes, $-a^2$ and $1$, and two non-real zeroes $\pm ia$, so $Z_{C}(p)=2$. The
function~$H_{\tfrac{3}{4}}[p]$ of this polynomial satisfies the following equality.

$$\frac{4}{3}H_{\tfrac34}[p](x)=(a^2-1)^2x^4-8a^2(a^2-1)x^3-2a^2(a^4-10a^2+1)x^2$$
$$ + 8a^4(a^2-1)x+a^4(a^2-1)^2$$
\begin{equation}
= \left[(a-1)x-a(a+1)\right]^2\cdot \left[(a+1)x+a(a-1)\right]^2.
\end{equation}
So, whenever $a\in\mathbb{R}\backslash\{-1,1\}$ it has four real  zeroes. Thus,
\begin{equation}\label{counter1.ineq}
Z_{\mathbb{R}}\left(H_{\tfrac{3}{4}}[p]\right)=4>2=Z_{\mathbb{C}}(p),
\end{equation}
and Conjecture~\ref{Shapiro1} fails.

Now we will prove Propositions~\ref{counterAnna} which gives a counterexample 
to Conjecture~\ref{Shapiro1} for each $n=\deg p \geq 5.$

{\bf Proof of Proposition~\ref{counterAnna}.} Let $n\geqslant5$. Consider the polynomial
\begin{equation}
p(x)=(x-1)^n+(x+1)^n.
\end{equation}
The numbers $z_k=i \cot\left(\frac{\pi}{2n}+\frac{\pi k}{n}\right), \   k=0,1,\ldots, n-1,$ are 
the zeros of this polynomial.  So, it has exactly $2\left\lfloor\dfrac n2\right\rfloor$ pure imaginary 
zeroes.

For the polynomial  $H_{\tfrac{n-1}{n}}[p]$, we have
\begin{equation}
H_{\tfrac{n-1}{n}}[p](x)=- 4n(n-1)(x-1)^{n-2}(x+1)^{n-2}.
\end{equation}
So, we have
\begin{equation}
Z_{\mathbb{R}}\left(H_{\tfrac{n-1}{n}}[p]\right)=2n-4>2
\left\lfloor\dfrac n2\right\rfloor=Z_{\mathbb{C}}(p)
\end{equation}
for any $n\geqslant5$. Thus, we proved that Conjecture~\ref{Shapiro1} is not true for 
any $n\geqslant4$.

At the end of this Section we will prove that the other Boris~Shapiro's conjecture~\ref{Shapiro2} 
is not valid. Let us prove Proposition~\ref{counterOlga2}.

{\bf Proof of Proposition~\ref{counterOlga2}.} Let $n\geq 2,$ and
$$p_{2n}(x)=\frac{x^{2n}}{2n}+\frac{x^2}{2}+1.$$
 It is obvious that $Z_{\mathbb{R}}(p)=0.$ For the polynomial $p_{2n}$ we have
\begin{equation}\label{counterO1}
-2n H_{\tfrac{2n-1}{2n}}[p](x) =x^{2n}(2n^2-5n+3)+2n(2n-1)x^{2n-2}-(n-1)x^2+2n.
\end{equation}

Let us show that the polynomial (\ref{counterO1}) doesn't have real zeros. First, assume that 
$|x|\geq 1.$ Since $n\geq 2,$ in this case $x^{2n-2}\geq x^2,$ and the following estimation is true
\begin{equation}\label{counterO2}
-2n H_{\tfrac{2n-1}{2n}}[p](x)\geq  (n-1)(2n-3)x^{2n} +(2n(2n-1)-(n-1))x^2+2n>0.
\end{equation}
If $|x|<1,$ then $x^2<1,$ and
\begin{equation}\label{counterO3}
-2n H_{\tfrac{2n-1}{2n}}[p](x)\geq   (n-1)(2n-3)x^{2n} +2n(2n-1)x^{2n-2}+n+1 >0.
\end{equation}
It follows from (\ref{counterO2}) and (\ref{counterO3}) that 
$$Z_{\mathbb{R}}\left(H_{\tfrac{2n-1}{2n}}[p](x)\right)=0.$$
Therefore,
$$Z_{\mathbb{R}}\left(H_{\tfrac{2n-1}{2n}}[p](x)\right)+Z_{\mathbb{R}}(p)=0,$$
and Boris~Shapiro's conjecture~\ref{Shapiro2} is not valid.

\section{ The function $M[p](x).$ Proof of Theorem~\ref{Theorem:lower.bound} and 
Theorem~\ref{conjectureSh2}}  \label{Sec4}

Without loss of generality we assume that $p$ is a monic polynomial, that is
\begin{equation}
\label{p1}
p(x)=x^n+a_{n-1}x^{n-1}+\ldots
\end{equation}

Let us introduce the following function
\begin{equation}
\label{ttss27}
M[p](x)=\frac{ p(x)p''(x)}{(p'(x))^2}.
\end{equation}

By virtue of (\ref{main.function.2})
\begin{equation}
\label{ttss28}
Q_{\varkappa}[p](x)=\varkappa-M[p](x).
\end{equation}
It is clear that  $Z_{\mathbb{R}}\left(Q_{\varkappa}[p]\right)$  is the 
same as the number of real roots of the equation
\begin{equation}
\label{mainequation}
M[p](x)=\varkappa.
\end{equation} 

Note that if $p'(x)$ doesn't have real roots, all the statements of 
Theorem~\ref{Theorem:lower.bound} 
and Theorem~\ref{conjectureSh2} are obviously true.
Let us denote by 
\begin{equation}\label{ends}
\xi_1<\xi_2<\cdots<\xi_{m-1}, \ \ m<n,
\end{equation}
different real poles of the function $M[p].$

\begin{remark}\label{multipleroots}
Note that if $\xi$ is a root of multiplicity $j$ of the polynomial $p,$ then
\begin{equation}\label{limit}
\lim_{x\to \xi} M[p](x)=\lim_{x\to \xi} \frac{ p(x)p''(x)}{(p'(x))^2}=\frac{j-1}{j}.
\end{equation}
So, the poles $\{\xi_k\}_ {k=1}^{m-1}$ from \ref{ends} are those zeros of $p'$ that are not roots of $p.$
\end{remark}

By virtue of Rolle's theorem each of the intervals below contains at most one zero of the polynomial $p$ 
\begin{equation}
\label{i1}
 I_1=(-\infty,\xi_1 ), \ I_2=(\xi_1,\xi_2 ),\ldots, \ I_m=(\xi_{m-1},\infty). 
\end{equation}

Let us divide all the intervals  $ (I_k)_{k=1}^{m}$ into two groups.

\begin{definition}\label{Intervals}
 {\it We will say that an interval }
\begin{itemize}
\item $I_k, \ k=1, 2, \ldots, m,$ {\it is an interval of the \textbf{first type} if it contains 
a root of the polynomial $p$;}
\item $I_k, \ k=1,2,\ldots, m,$ {\it is an interval of the \textbf{second type} 
if it doesn't contain  roots of the polynomial $p.$}
\end{itemize}
\end{definition}

Assume that $\xi$ is a root of the derivative $p'$ which is not a root of $p.$ Then
\begin{equation}\label{expanp}
p(x)=p(\xi)+\frac{p^{(s)}(\xi)}{s!}(x-\xi)^s+\ldots, \ \ \mbox{where} \ \ p^{(s)}(\xi) \neq 0.
\end{equation}

If we take up the position that the right behavior of a polynomial $p$ is to have a 
root in the interval $I_k,$ while not to have a root in the interval is a wrong behavior, 
we will come to the following classification of  the points $\{\xi_k \}_{k=1}^{m-1}$ 
(see \ref{expanp}).

\begin{definition}\label{endpoints}
 Assume that $\xi$ is a pole of $M[p]$ (see (\ref{ttss27})).
\begin{itemize}
\item {\it $\xi$ is a right-hand~\textbf{right point} if 
$$p(\xi) p^{(s)}(\xi)(x-\xi)^s<0, \ \mbox{when} \ x \ \mbox{is close to}
\ \xi \ \mbox{and} \  x>\xi.$$
Thus, } 
\begin{equation}\label{IRR}
\lim_{x \to \xi^{+}} M[p](x)=-\infty.
\end{equation}
\item {\it $\xi$ is a left-hand \textbf{right point} if 
$$p(\xi) p^{(s)}(\xi)(x-\xi)^s<0, \ \mbox{when} \ x \ \mbox{is close to}
\ \xi \ \mbox{and} \  x<\xi.$$
 Thus,  }
\begin{equation}\label{ILR}
\lim_{x \to \xi^{-}} M[p](x)=-\infty.
\end{equation}
\item  {\it $\xi$ is a right-hand \textbf{wrong point} if
 $$p(\xi) p^{(s)}(\xi)(x-\xi)^s>0, \ \mbox{when} \ x \ \mbox{is close to}\ \xi \ \mbox{and} \  x>\xi.$$
 Thus,  }
\begin{equation}\label{IRW}
\lim_{x \to \xi^{+}} M[p](x)=+\infty.
\end{equation}
\item {\it $\xi$ is a left-hand \textbf{wrong point} if\\ $$p(\xi) p^{(s)}(\xi)(x-\xi)^s>0, \ 
\mbox{when} \ x \ \mbox{is close to}\ \xi \ \mbox{and} \  x<\xi.$$
 Thus,  }
\begin{equation}\label{ILW}
\lim_{x \to \xi^{-}} M[p](x)=+\infty.
\end{equation}
\end{itemize} 
\end{definition}

\begin{remark}\label{Intend}
The above definition says the following about the end-points of the intervals  $I_k.$ 
\begin{itemize}
\item {\it If a finite interval $I_k=(\xi_{k-1}, \xi_k)$ is an 
interval of the first type, then $\xi_{k-1}$ is a right-hand right point, and $\xi_k$ is 
a left-hand right point.}
\item {\it If a finite interval $I_k=(\xi_{k-1}, \xi_k)$ is an 
interval of the second type, then either $\xi_{k-1}$ is a right-hand right point, and 
$\xi_k$ is a left-hand wrong point, or $\xi_{k-1}$ is a right-hand wrong point, and 
$\xi_k$ is a left-hand right point.}
\item {\it If an infinite interval $I_1=(-\infty, \xi_1)$ is an interval of the first type, 
its finite endpoint $\xi_1$ is a left-hand right point. If an infinite interval 
$I_m=(\xi_{m-1}, \infty)$ is an interval of the first type, its finite endpoint 
$\xi_{m-1}$ is a right-hand right point. }
\item {\it If an infinite interval $I_1=(-\infty, \xi_1)$ is an interval of the second 
type, its finite endpoint $\xi_1$ is a left-hand wrong point. If an infinite interval 
$I_m=(\xi_{m-1}, \infty)$ is an interval of the second type, its finite endpoint 
$\xi_{m-1}$ is a right-hand wrong point.}
\end{itemize}
\end{remark}

So, using (\ref{ttss28}), (\ref{mainequation}), (\ref{IRR})-(\ref{ILW}) we come to the 
following conclusion about the number of real roots of the function $Q_{\varkappa}[p]$ 
on finite intervals $I_k.$

\begin{lemma}\label{Interval.endpoints}
Let $p$ be a real polynomial of degree $n\geq 2$. Then
\begin{itemize}
\item  If a finite interval $I_k$ is an interval of the first type, 
then for any real $\varkappa$ the number $Z_{I_k}\left(Q_{\varkappa}[p]\right)$ 
is even.
\item If a finite interval $I_k$ is an interval of the second type, 
then for any real $\varkappa$ the number $Z_{I_k}\left(Q_{\varkappa}[p]\right)$ 
is odd.
\end{itemize}
\end{lemma}

Let us discuss the behavior of the function $M[p]$ at infinity. Note that by (\ref{ttss27}) we have
\begin{equation}\label{infinity}
\lim_{x\to \pm \infty}M[p](x)=\frac{n-1}{n}.
\end{equation}

The statement below follows from (\ref{infinity}), two last statements of Remark~\ref{Intend} 
and  (\ref{IRR})-(\ref{ILW}). It gives an estimation for  the number of real roots of the function 
$Q_{\varkappa}[p]$ on infinite intervals $I_k.$ As we noted in the beginning of this section, 
the case when $p'$ doesn't have real roots is trivial for Theorem~\ref{Theorem:lower.bound} 
and Theorem~\ref{conjectureSh2}.  That's why we do not consider it in the lemma below.

\begin{lemma}\label{Infint}
Let $p$ be a real polynomial of degree $n\geq 2$, and its derivative $p'$ has at least one 
real root.  Then 
\begin{itemize}
\item  If an infinite interval $I_k, \  k=1 \  \mbox{or}\  m,$  is an interval of the first type, 
and $\varkappa \geq \frac{n-1}{n},$ then the number $Z_{I_k}\left(Q_{\varkappa}[p]\right)$ 
is even.
\item  If an infinite interval $I_k,  \ k=1 \  \mbox{or} \ m,$  is an interval of the first type, and 
$\varkappa <  \frac{n-1}{n},$ then the number $Z_{I_k}\left(Q_{\varkappa}[p]\right)$ is odd.
\item  If an infinite interval $I_k, \ k=1 \  \mbox{or} \  m,$  is an interval of the second type, and 
$\varkappa > \frac{n-1}{n},$ then the number $Z_{I_k}\left(Q_{\varkappa}[p]\right)$ is odd.
\item  If an infinite interval $I_k, \ k=1 \  \mbox{or}\ m,$  is an interval of the second type, and 
 $\varkappa \leq \frac{n-1}{n},$ then the number $Z_{I_k}\left(Q_{\varkappa}[p]\right)=
 Z_{I_k}\left(H_{\varkappa}[p]\right)$ is even.
\end{itemize}
\end{lemma}

Now we can prove Theorem~\ref{Theorem:lower.bound} on lower bounds for the number of 
real roots of $Q_{\varkappa}[p].$ We will use the following standard notation.
\begin{notation}\label{numberelements}
Let $S$ be a finite set. Then by $\#S$ we will denote the number of elements in the set $S$.
\end{notation}

{\bf Proof of Theorem~\ref{Theorem:lower.bound}.} Let $\varkappa > \frac{n-1}{n}.$ We 
need to estimate from below $Z_{\mathbb{R}}\left(Q_{\varkappa}[p]\right)$. 
By virtue of Lemma~\ref{Interval.endpoints} and Lemma~\ref{Infint}, if we consider a 
contribution of each $Z_{I_k}\left(Q_{\varkappa}[p]\right)$ into $Z_{\mathbb{R}}
\left(Q_{\varkappa}[p]\right)$ as minimum possible, we obtain
$$Z_{\mathbb{R}}\left(Q_{\varkappa}[p]\right) \ \geq \ \#
 \{I_k| \  I_k \ \mbox{is an interval of the second type} \} \ $$ 
$$= m \ - \  \#\{ I_k| \ I_k  \mbox{ is an interval of the first type}\}. $$
Recall that by (\ref{ends}), (\ref{limit}), (\ref{i1}) the number of pairwise unequal roots of 
the polynomial $p'$ that are not roots of $p $ equals $m-1.$ So,
\begin{equation}\label{equ1}
 m \ = \ \mathring{Z}_\mathbb{R}(p') \ + 1 \ - \ \# \{\xi| \ \xi \ 
 \mbox{is a multiple root of $p$}\}.
\end{equation}
Note that
$$\#\{ I_k| \ I_k  \mbox{ is an interval of the first type}\} \ =$$
\begin{equation}\label{equ2}
 \mathring{ Z}_\mathbb{R}(p) \ - \ \# \{\xi| \ \xi \ \mbox{is a multiple root of $p$}\}.
\end{equation}
The first statement of Theorem~\ref{Theorem:lower.bound} follows immediately from 
(\ref{equ1}) and (\ref{equ2}).

Now the second statement of Theorem~\ref{Theorem:lower.bound} is obvious. Counting the 
number of the intervals of the second type we need just to omit two infinite intervals.

In order to prove the third statement of Theorem~\ref{Theorem:lower.bound} we need 
the following lemma.
\begin{lemma}\label{IntI1} Let $p$ be a real polynomial of degree $n\geq 2$. If $I_k$ 
is a finite interval of the first type, and $\varkappa \leq 0,$ then
\begin{equation}\label{IntervalI1}
Z_{I_k}\left(Q_{\varkappa }[p](x)\right)\geq 2.
\end{equation}
\end{lemma}

{\bf Proof of Lemma~\ref{IntI1}.} By virtue of Rolle's theorem inside each 
finite interval $I_k$ there is a root $\gamma$ of $p''.$ Thus, by (\ref{ttss27}) if 
$I_k$ is an interval of the first type, it contains at least two roots of the function $M[p].$ 
Therefore, taking into account formulas~(\ref{IRR}), (\ref{ILR}) and (\ref{ttss28}), we obtain
$$Z_{I_k}\left(\varkappa -M[p](x)\right)=Z_{I_k}\left(Q_{\varkappa }[p](x)\right)\geq 2.$$
Lemma~\ref{IntI1} is proved. $\square$

It follows from Lemma~\ref{Interval.endpoints}, Lemma~\ref{Infint} and Lemma~\ref{IntI1} that
$$Z_{\mathbb{R}}\left(Q_{\varkappa }[p](x)\right)\geq
\#\{I_k | \ I_k \ \mbox{is a finite interval of the second type}\} \ +$$
$$+\ 2\cdot \# \{I_k| \ I_k \ \mbox{is a finite interval of the first type}\}
\ + $$  $$\#\{I_k| \ I_k \ \mbox{is an infinite interval of the first type}\} \ =$$
$$ \# \{I_k| \ I_k \ \mbox{is a finite interval}\}
+\#\{I_k| \ I_k \ \mbox{is an interval of the first type}\}=$$  $$m-2+Z_{\mathbb{R}}(p).$$

The third statement of Theorem~\ref{Theorem:lower.bound} follows from (\ref{equ1}) 
and (\ref{equ2}). $\square$

{\bf Proof of Theorem~\ref{conjectureSh2}.} Note that if a polynomial $p$ has real 
roots, the statement is obvious.

If $\xi$ is a multiple real root of the polynomial $p',$ then
$$p(x)=a_0+a_m(x-\xi)^m q(x), \ \ m\geq 3,$$
and whence $\xi$ is a root of $H_{\frac{n-1}{n}}[p](x)=\frac{n-1}{n} (p'(x))^2 - p(x)p''(x).$ 
So, in this case theorem is also obvious.

Assume now that $p$ doesn't have real roots, all real roots of $p'$ are simple, and
$p(x)=x^n+a x^{n-1}+\ldots$
Let us expand $p$ about $-\frac{a}{n}.$
\begin{equation}
\label{abd7}
 p(x)=x^n+a x^{n-1}+\ldots=\left(x+\frac{a}{n}\right)^n+
 \ b\left(x+\frac{a}{n}\right)^{n-2}+\ldots
\end{equation} 
Thus,
\begin{equation}
\label{abd8}
p'(x)=n \left(x+\frac{a}{n}\right)^{n-1}+\ b(n-2)\left(x+\frac{a}{n}\right)^{n-3}+\ldots
\end{equation} 
 and  
\begin{equation}
\label{abd9}
p''(x)=n(n-1) \left(x+\frac{a}{n}\right)^{n-2}+\ b(n-2)(n-3)\left(x+\frac{a}{n}\right)^{n-4}+\ldots
\end{equation}
It follows from (\ref{abd7}), (\ref{abd8}) and (\ref{abd9}) that
\begin{equation}
\label{Identity}
M[p](x)-\frac{n-1}{n}=\frac{2b n  \left(x+\frac{a}{n}\right)^{2n-4}+\ldots}{\left(n 
\left(x+\frac{a}{n}\right)^{n-1}+b(n-2)\left(x+\frac{a}{n}\right)^{n-3}+\ldots\right)^2}.
\end{equation}
It is given that for some $k=1,2,\ldots, n-2,$ all the zeros of the derivative
 $$p^{(k)}(x)=\frac{n!}{(n-k)!}\left(x+\frac{a}{n}\right)^{n-k}+b\  
 \frac{(n-2)!}{(n-k-2)!}\left(x+\frac{a}{n}\right)^{n-k-2}+\cdots$$ 
are real. Therefore, the coefficient $b$ in (\ref{abd7}) must be negative 
(see~\cite[p. 434, Lemma 3]{Levin}). Since $b<0,$  by (\ref{Identity}) 
there exists such a value $x_0$ that for $|x|>x_0>0$ the following is true
\begin{equation}
\label{abd10}
M[p](x)-\frac{n-1}{n}<0.
\end{equation}
It means that the graph of the function $y=M[p](x)$ approaches its asymptote $y=\frac{n-1}{n}$ from below. 

Therefore, the statement of Theorem~\ref{conjectureSh2} follows immediately from the last statement of 
Remark~\ref{Intend} and formulas~(\ref{IRW}), (\ref{ILW}). $\square$

\section{ The functions $\varphi_{\lambda}(x)$ and $\Phi_{\lambda}(x)$ } \label{Sec5}

From now on we will assume that all real roots of the polynomial $p$ are simple, and all roots of 
its derivative $p'$ are real and simple. So, in this case the number of intervals $I_k$ is equal to $n.$ 
The following six sections are devoted to the proof of Theorem~\ref{Hawaiitrue}. First we are 
going to study the case of $\varkappa>\frac{n-1}{n}.$  Put
\begin{equation}
\label{b2}
 \lambda=\varkappa -1.
\end{equation}
Consider the function
\begin{equation}
\label{b1}
 \varphi_{\lambda} (x)=\frac{p(x)}{p' (x) }+\lambda x.
\end{equation} 
 Note that 
\begin{equation}
\label{b3}
\varphi_{\lambda}'(x)=\frac{(p' (x))^2-p(x) p'' (x)}{(p' (x))^2 }+\lambda=
\frac{\varkappa (p' (x))^2-p(x) p'' (x)}{(p'(x))^2 }=Q_{\varkappa}[p](x).
\end{equation}
Therefore,
\begin{equation}
\label{a1}
 Z_{\mathbb{R}}(\varphi_{\lambda}'(x))=Z_{\mathbb{R}} \left(Q_{\varkappa}[p](x)\right).
\end{equation}
If $(\xi_k, p(\xi_k))$ is a local maximum, then $p'$ changes sign from positive to negative, 
and if $(\xi_k, p(\xi_k))$ is a local minimum, then $p'$ changes sign from negative to positive as 
$x$ increases through $\xi_k$. Therefore, according to definition~\ref{endpoints} and 
(\ref{b1}), if $\xi_k$ is a right point, then
\begin{equation}
\label{a2}
\lim_{x \to \xi_k^{+}}\varphi_{\lambda} (x)=-\infty \ \ \mbox{and} \ \ \lim_{x \to 
\xi_k^{-}}\varphi_{\lambda} (x)=+\infty,
\end{equation}
 and if $\xi_k$ is a wrong point, then
\begin{equation}
\label{a3}
\lim_{x \to \xi_k^{+}}\varphi_{\lambda} (x)=+\infty \ \ \mbox{and} \ \ 
\lim_{x \to \xi_k^{-}}\varphi_{\lambda} (x)=-\infty. 
\end{equation}

The formula below describes asymptotic behavior of the function 
$ \varphi_{\lambda} (x) $ at infinity
\begin{equation}
\label{c1}
\varphi_{\lambda} (x)=\frac{p(x)}{p' (x) }+\lambda x=\left( \frac{1}{n}+\lambda\right)x+
o(x), \ \ x \to \pm \infty. 
\end{equation}
 Given a real number $\lambda,$ we will evaluate the number of 
 all real roots of the equation 
\begin{equation}
\label{me1}
\varphi_{\lambda} (x)=\frac{p(x)}{p' (x) }+\lambda x=\mu
\end{equation}
in each interval $I_k$ for any real value of the 
parameter $\mu.$

The lemma below follows immediately from (\ref{a2}), (\ref{a3}), (\ref{c1}) and 
Remark~\ref{Intend}.
\begin{lemma}\label{lemma1} The following statements are true.
\begin{enumerate}
\item {\it For all real values of the parameters $\mu, \lambda,$ and every $k=2, 3, \ldots, n-1$}
\begin{itemize}
\item {\it $Z_{I_k}(\varphi_{\lambda}(x)-\mu)$ is odd if $I_k$ is an interval of the first type.}
\item {\it $Z_{I_k}(\varphi_{\lambda}(x)-\mu)$ is even if $I_k$ is an interval of the second type.}
\end{itemize}

\item {\it For all real values of the parameter $\mu$, for every $\lambda>-\frac{1}{n},$ 
and $k=1$ or $k=n$}
\begin{itemize}
\item {\it $Z_{I_k}(\varphi_{\lambda}(x)-\mu)$ is odd if $I_k$ is an interval of the first type.}
\item {\it $Z_{I_k}(\varphi_{\lambda}(x)-\mu)$ is even if $I_k$ is an interval of the second type.}
\end{itemize}

\item  {\it For all real values of the parameter $\mu$, for every $\lambda<-\frac{1}{n},$ 
and $k=1$ or $k=n$}
\begin{itemize}
\item {\it $Z_{I_k}(\varphi_{\lambda}(x)-\mu)$ is even if $I_k$ is an interval of the first type.}
\end{itemize}

\end{enumerate}
\end{lemma}

Further we will need the following obvious lemma that describes some properties of a monic 
polynomial $p,$ all roots of whose derivative $p'$ are real and simple. 
\begin{lemma} \label{lemma2} Let $p$ be a real monic polynomial, and all zeros of $p'$ are real 
and simple. Then
\begin{enumerate}
\item {\it $(-1)^{n-k-1} p''(\xi_k)>0, \ \  k=1,2,\ldots, n-1.$}
\item  {\it $(-1)^{n-k} p'(x)>0,  \ \  x\in I_k.$ }
\item {$\lim_{x \to \infty} p'(x)=\lim_{x \to -\infty} (-1)^{n-1} p'(x)=+\infty.$}
\end{enumerate}
\end{lemma}

In order to estimate the number of real roots of (\ref{me1}) in intervals $I_k$ 
we will consider a function
\begin{equation}
\label{aa9}
 \Phi_{\lambda}(x, \mu)=(\varphi_{\lambda} (x)-\mu)p' (x)  =p(x)+\lambda x p'(x)-\mu p'(x),
\end{equation} 
and its derivative
\begin{equation}
\label{aa10}
 \Phi_{\lambda}'(x, \mu)=(1+\lambda) p'(x)+(\lambda x -\mu) p''(x).
\end{equation} 

\section{ Proof of Theorem~\ref{Hawaiitrue} in the Partial Case of  $ \lambda=0$ } \label{Sec6}

Let us prove Theorem~\ref{Hawaiitrue} in the partial case of $\lambda=0.$ In a sense, this is the 
base case in the proof of this theorem. 

If $\lambda=0, $ then (\ref{b1}) has the form
\begin{equation}
\label{le1}
\varphi_{0} (x)=\frac{p(x)}{p' (x) }.
\end{equation} 
The statement below describes monotonicity of the function $\varphi_{0},$ provides estimations 
for the number of real roots of the derivative $\varphi'_{0}$ on the intervals $I_k,$  and as a 
result gives the proof of Theorem~\ref{Hawaiitrue} in the case of $\lambda=0.$

\begin{lemma}\label{phi0}
The following statements are true.
\begin{enumerate}
\item  The function $ \varphi_{0} (x)=\frac{p(x)}{p' (x) }$  is increasing on any interval $I_k$ 
of the first type.
\item  $Z_{I_k}(\varphi_{0}'(x))= Z_{I_k}((p')^2-p p'')=0$ if $I_k$ is an interval of the first 
type, and 
\item $Z_{I_k}(\varphi_{0}'(x))=Z_{I_k}((p')^2-p p'')=1$ if $I_k$ is an interval of the 
second type.
\item  $Z_{\mathbb{R}}(\varphi_{0}'(x))= Z_{\mathbb{R}}((p')^2-p p'')=Z_{\mathbb{C}}(p).$ 
\end{enumerate}
\end{lemma}

{\bf Proof of Lemma~\ref{phi0}.} If $\lambda =0,$ then (\ref{aa9}) and  (\ref{aa10}) can be 
rewritten in the form
\begin{equation}
\label{le2}
 \Phi_{0}(x, \mu)=(\varphi_{0} (x)-\mu)p' (x) =p(x)-\mu p'(x),
\end{equation}
and
\begin{equation}
\label{le3}
\Phi_{0}'(x, \mu)= p'(x) -\mu p''(x). 
\end{equation} 
Since $\xi_k$ is a root of the polynomial $p'$ the following identity is true
\begin{equation}
\label{aa11}
(-1)^{n-k} \Phi_{0}'(\xi_k, \mu)=( -1)^{n-k+1} \mu p''(\xi_k), \ \ k=1,2,\ldots, n-1.
\end{equation}
It follows from the statement 1 of Lemma~\ref{lemma2} that
\begin{equation}
\label{aa12}
(-1)^{n-k} \Phi_{0}'(\xi_k, \mu)>0, \ \ \mbox{whenever} \ \ \mu>0,
\end{equation} 
and
\begin{equation}
\label{aa13}
(-1)^{n-k} \Phi_{0}'(\xi_k, \mu)<0, \ \ \mbox{whenever} \ \ \mu<0.
\end{equation} 
By the statement 3 of Lemma~\ref{lemma2}
\begin{equation}
\label{aa14}
\lim_{x \to \infty}\Phi_{0}'(x, \mu)=\lim_{x \to -\infty} (-1)^{n-1}\Phi_{0}'(x, \mu)=+\infty.
\end{equation}
Since $\deg \Phi_{0}'(x, \mu)=n-1,$ it follows from (\ref{aa12}) and (\ref{aa14}) that 
for any $\mu>0$ 
\begin{equation}
\label{aa15}
Z_{I_k}(\Phi_{0}'(x, \mu))=1, \ \ \mbox{when} \ \ k=2, 3, \ldots, n,
\end{equation}
and from (\ref{aa13}) and (\ref{aa14}) that for any $\mu<0$ 
\begin{equation}
\label{aa16}
Z_{I_k}(\Phi_{0}'(x, \mu))=1, \ \ \mbox{when} \ \ k=1, 2, \ldots, n-1.
\end{equation}
Thus, for any real $\mu\neq 0$ and $k=1,2,\ldots,n,$ we obtain
\begin{equation}
\label{aa17}
Z_{I_k}(\Phi_{0}(x, \mu))=Z_{I_k}(\varphi_{0}(x)- \mu) \leq 2. 
\end{equation}
If $\mu=0,$ then according to $(\ref{le2})$
$$\Phi_{0}(x, 0)=\varphi_0(x)=p(x).$$ 
Hence, by virtue of Rolle's Theorem
$$Z_{I_k}(\Phi_{0}(x, 0))=Z_{I_k}(\varphi_{0}(x))=Z_{I_k}(p(x))  \leq 1, $$
and (\ref{aa17}) is obviously true.

It follows from (\ref{aa17}), and the statement 1 of Lemma~\ref{lemma1} that for any 
real $\mu$ and each interval $I_k$ of the first type we have
\begin{equation}
\label{aa18}
Z_{I_k}(\Phi_{0}(x, \mu))=Z_{I_k}(\varphi_{0}(x)- \mu) =1. 
\end{equation}

According to the Remark~\ref{Intend}, both endpoints of a finite interval of the first type 
$I_k$ are right points;  and the finite endpoint of an infinite interval of the first type $I_k$ is 
a right point. So, if $I_k$ is an interval of the first type,  the relation  (\ref{aa18}) together with 
(\ref{a2})  in the case of a finite interval $I_k,$ and (\ref{a2}),  (\ref{c1}) in the case of an 
infinite interval $I_k,$ provide the fact that the function $\varphi_{0}(x)$ is increasing on each 
interval $I_k$ of the first type. The statement~1 of Lemma~\ref{phi0} is proved.

Since $\varphi_{0}(x)$ is monotone on each interval $I_k$ of the first type, we have
\begin{equation}
\label{aa19}
Z_{I_k}(\varphi_{0}'(x))=Z_{I_k}((p'(x))^2-p(x)p''(x)) =0,
\end{equation}
where $I_k$ is an interval of the first type.
It follows from the statement 2 of Lemma~\ref{lemma1} and (\ref{aa17}) that for any 
real $\mu$ and each interval $I_k$ of the second type
\begin{equation}
\label{aa20}
Z_{I_k}(\Phi_{0}(x, \mu))=Z_{I_k}(\varphi_{0}(x) - \mu) =0  \ \mbox{or} \ 2.
\end{equation}
Using the statements of the Remark~\ref{Intend} related to intervals of the second type, 
taking into account (\ref{a2}) and (\ref{a3}) in the case of a finite interval $I_k,$ and  
(\ref{a3})(\ref{c1}) in the case of an infinite interval $I_k,$ we obtain that for each 
interval $I_k$ of the second type
\begin{equation}
\label{aa21}
Z_{I_k}(\varphi_{0}')=Z_{I_k}((p')^2-p  p'') =1.
\end{equation}
The statement 2 of Lemma~\ref{phi0} is proved.

The formulas (\ref{aa19}) and (\ref{aa21}) show that the number of real zeros of the 
polynomial $(p')^2-pp''$ coincides with the total number of the intervals of the 
second type, that is 
\begin{equation}
\label{aa22}
Z_{\mathbb{R}}((p')^2-p  p'') =n-Z_{\mathbb{R}}(p)=Z_{\mathbb{C}}(p).
\end{equation}
Lemma~\ref{phi0} is proved. $\square$

\section{The function $\Psi_{\alpha}(x, \nu)$} \label{Sec7}

Now we assume that $\lambda \neq 0.$ So, we can rewrite (\ref{aa10}) in the form 
\begin{equation}
\label{aa23}
 \frac{1}{\lambda}\Phi_{\lambda}'(x, \mu)=\frac{1+\lambda}{\lambda} p'(x)+
 \left( x -\frac{\mu}{\lambda}\right) p''(x),
\end{equation}
where $\mu$ is any real number.

Denote by 
\begin{equation}
\label{aa24}
 \alpha=\frac{1+\lambda}{\lambda}, 
\end{equation} 
 and by
\begin{equation}
\label{aaa24}
\nu=\frac{\mu}{\lambda}. 
\end{equation} 
Note that $\nu$ takes any real value, while
\begin{equation}
\label{aa25}
 \alpha>1,  \ \ \mbox{when} \ \ \lambda>0,
\end{equation} 
and 
\begin{equation}
\label{aa26}
 \alpha<-(n-1),  \ \ \mbox{when} \ \ 1/n<\lambda<0.
\end{equation} 

We introduce the function
\begin{equation}
\label{aa27}
 \Psi_{\alpha}(x, \nu)=:\alpha p'(x)+\left( x -\nu \right) p''(x).
\end{equation} 
Obviously,
\begin{equation}
\label{aa28}
\Psi_{\frac{1+\lambda}{\lambda}}\left(x, \frac{\mu}{\lambda}\right) =
\frac{1}{\lambda}\Phi_{\lambda}'(x, \mu). 
\end{equation}

The lemma below describes behavior of the function $\Psi_{\alpha}(x, \nu)$ for several 
important values of $x$ and at infinity. 

\begin{lemma}\label{lemma4} The  following statements are true.
\begin{enumerate}
\item {\it For all real values of the parameter $\nu$ and every $k=1, 2, \ldots, n-1$}
\begin{itemize}
\item {\it $(-1)^{n-k}\Psi_{\alpha}(\xi_k, \nu)<0,$ whenever $\xi_k>\nu.$} 
\item {\it $(-1)^{n-k}\Psi_{\alpha}(\xi_k, \nu)>0,$ whenever $\xi_k<\nu.$}
\end{itemize}
\item {\it If $\nu \in I_k=(\xi_{k-1}, \xi_k), \ \ k=2, \ldots, n-1$, then}
\begin{itemize}
\item {\it  $(-1)^{n-k}\Psi_{\alpha}(\nu, \nu)>0,$ whenever $\alpha>0.$}
\item {\it  $(-1)^{n-k}\Psi_{\alpha}(\nu, \nu)<0,$ whenever $\alpha<0.$}
\end{itemize}
\item  {\it For all real values of the parameter $\nu$} 
\begin{itemize}
\item {\it $\lim_{x\to +\infty}\Psi_{\alpha}(x, \nu)=\lim_{x\to -\infty} (-1)^{n-1}\Psi_{\alpha}(x, \nu)=
+\infty,$ whenever $\alpha >-(n-1)$.}
\item {\it $\lim_{x\to +\infty}\Psi_{\alpha}(x, \nu)=\lim_{x\to -\infty} (-1)^{n-1}
\Psi_{\alpha}(x, \nu)=-\infty,$ whenever $\alpha <-(n-1)$.}
\end{itemize}
\item {\it  For all real values of the parameter $\alpha$ and any $k=1, 2, \ldots, n-1$}
$$\Psi_{\alpha}(x, \xi_k)=(x-\xi_k)h(x), $$
{\it where}
$$h(\xi_k)=(\alpha +1) p'' (\xi_k).$$
\end{enumerate}
\end{lemma}

{\bf Proof of Lemma~\ref{lemma4}.}  Since $\xi_k$ is a root of $p',$ by virtue 
of (\ref{aa27}) we have

\begin{equation}
\label{aa29}
 (-1)^{n-k}\Psi_{\alpha}(\xi_k, \nu)= (-1)^{n-k-1} p''(\xi_k) (\nu-\xi_k).
\end{equation} 
Now the first statement of Lemma~\ref{lemma4} follows immediately from the first 
statement of Lemma~\ref{lemma2}.

The following identity can also be derived from (\ref{aa27})
\begin{equation}
\label{aa36}
(-1)^{n-k} \Psi_{\alpha}(\nu, \nu)=(-1)^{n-k}\alpha p'(\nu) .
\end{equation} 
The second statement of Lemma~\ref{lemma4} is a direct corollary of this identity and 
the second statement of Lemma~\ref{lemma2}.

Using definition (\ref{aa27}) of the function $\Psi_{\alpha}$ and the third statement of 
Lemma~\ref{lemma2} we easily obtain the third statement of Lemma~\ref{lemma4}.

Let us prove the fourth statement of Lemma~\ref{lemma4}. Since $\xi_k$ is a root of the polynomial $p',$ 
there exists a polynomial $g$ such that
\begin{equation}
\label{aa32}
p'(x)=(x-\xi_k)g(x),
\end{equation} 
and 
\begin{equation}
\label{aa33}
g(\xi_k)=\lim_{x \to \xi_k}\frac{p'(x)-p'(\xi_k)}{x-\xi_k} =p''(\xi_k).
\end{equation}
It follows from (\ref{aa27}) with $\nu=\xi_k$ and (\ref{aa32}) that
\begin{equation}
\label{aa34}
\Psi_{\alpha}(x, \xi_k)=\alpha p'(x)+(x-\xi_k) p''(x)= (x-\xi_k)(\alpha g(x)+p''(x))=:(x-\xi_k)h(x),
\end{equation} 
where
\begin{equation}
\label{aa35}
h(\xi_k)=(\alpha +1)p''(\xi_k).
\end{equation}

Lemma~\ref{lemma4} is proved. $\square$

\section{Estimations for $Z_{I_k}\left(\Psi_{\alpha}(x, \nu)\right)$} \label{Sec8}

The following Lemma  provides the main technical tools for proving Theorem~\ref{Hawaiitrue}.
\begin{lemma}\label{lemma5} The following statements are true.
\begin{enumerate}
\item If $\alpha \geq 0,$ then for all real values of the parameter $\nu$ we have
\begin{itemize}
\item $Z_{I_k}(\Psi_{\alpha}(x, \nu))\leq 2$ when $k=2, 3, \ldots, n-1.$
\item  $Z_{I_1}(\Psi_{\alpha}(x, \nu))\leq 1$ and $Z_{I_n}(\Psi_{\alpha}(x, \nu))\leq 1.$
\end{itemize}
\item If $\alpha < -(n-1),$ then for all real values of the parameter $\nu$ and all 
$k=1, 2, \ldots, n$ we have   $Z_{I_k}(\Psi_{\alpha}(x, \nu))\leq 1.$
\end{enumerate}
\end{lemma}

{\bf Proof of Lemma~\ref{lemma5}.} First, consider the case of $\alpha=0.$ 
Using (\ref{aa27}) we obtain
$$\Psi_{\alpha}(x, \nu))=(x-\nu)p''(x),$$
and the statement 1 of Lemma~\ref{lemma5} follows from Rolle's theorem.

By definition (\ref{aa27}) of the function $\Psi_{\alpha}(x, \nu),$ for all real values of the 
parameter $\nu$ and $\alpha\neq -(n-1)$ we have
\begin{equation}
\label{aa37}
\deg\Psi_{\alpha}(x, \nu))=n-1. 
\end{equation} 

Now assume that $\alpha>0.$ Using the first statement of Lemma~\ref{lemma4} in the case of a 
finite interval $I_k,$ and additionally the third statement of Lemma~\ref{lemma4} in the case of 
the infinite intervals $I_1$ and $I_n,$ we conclude that the relations below are valid
\begin{equation}
\label{aa38}
 Z_{I_k}(\Psi_{\alpha}(x, \nu))=1 , \ \ k=2, 3, \ldots, n,\  \ \mbox{whenever} \ \  \nu>\xi_{n-1},
\end{equation} 
and 
\begin{equation}
\label{aaa38}
Z_{I_k}(\Psi_{\alpha}(x, \nu))=1,  \  \  k=1, 2, \ldots, n-1,\  \ \mbox{whenever} \ \ \nu<\xi_1.
\end{equation}

If $\xi_{j-1}<\nu<\xi_j,$ then applying the first and the second statements of 
Lemma~\ref{lemma4} we obtain the following relations.
\begin{equation}
\label{aa39}
Z_{I_k}(\Psi_{\alpha}(x, \nu))=1,  \  \  k=2,\ldots,j-1, j+1, \ldots, n-1,
\end{equation} 
and
\begin{equation}
\label{aaa39}
 Z_{(\xi_{j-1}, \nu)}(\Psi_{\alpha}(x, \nu))= Z_{( \nu, \xi_j)}(\Psi_{\alpha}(x, \nu))=1.
\end{equation} 
Therefore, 
\begin{equation}
\label{aaa40}
 Z_{I_1}(\Psi_{\alpha}(x, \nu))= Z_{I_n}(\Psi_{\alpha}(x, \nu))=0.
\end{equation}
This provides the first statement of Lemma~\ref{lemma5} in the case when 
$\nu\neq \xi_j, \ j=1, 2, \ldots, n-1.$

To investigate the case of $\nu=\xi_j, \ j=1, 2, \ldots, n-1,$ we will use the fourth statement 
of Lemma~\ref{lemma4}. According to this statement 
\begin{equation}
\label{bb29}
 \Psi_{\alpha}(x, \xi_j)=(x-\xi_j)h(x),
\end{equation}
where
\begin{equation}
\label{aa41}
h (\xi_j)=(\alpha+1)p''(\xi_j).
\end{equation} 
By the first statement of Lemma~\ref{lemma2}
\begin{equation}
\label{aa43}
(-1)^{n-j-1} h(\xi_j)=(-1)^{n-j-1}(\alpha+1)p''(\xi_j)>0, \ \ j=1,2, \ldots, n-1.
\end{equation}

It follows from (\ref{bb29}) and (\ref{aa43}) that there exists a positive number $\varepsilon$ 
such that
\begin{equation}
\label{aa44}
(-1)^{n-j-1} \Psi_{\alpha}(\xi_j+\varepsilon, \xi_j)=\varepsilon (-1)^{n-j-1} h(\xi_j) >0,
\end{equation} 
and
\begin{equation}
\label{aa45}
(-1)^{n-j-1} \Psi_{\alpha}(\xi_j-\varepsilon, \xi_j)=-\varepsilon (-1)^{n-j-1} h(\xi_j) <0 ,
\end{equation} 
while
\begin{equation}
\label{aaa45}
 \Psi_{\alpha}(\xi_j, \xi_j)=0.
\end{equation} 
Applying the first statement of Lemma~\ref{lemma4} we obtain
\begin{equation}
\label{aa46}
(-1)^{n-k} \Psi_{\alpha}(\xi_k, \xi_j)<0, \ \ \mbox{when} \ \ k > j,
\end{equation} 
and
\begin{equation}
\label{aa47}
(-1)^{n-k} \Psi_{\alpha}(\xi_k,  \xi_j)>0, \ \ \mbox{when} \ \ k < j.
\end{equation} 
In particular,
\begin{equation}
\label{aaa46}
(-1)^{n-j-1} \Psi_{\alpha}(\xi_{j+1}, \xi_j)<0 ,
\end{equation} 
and
\begin{equation}
\label{aaa47}
(-1)^{n-j-1} \Psi_{\alpha}(\xi_{j-1},  \xi_j)>0 .
\end{equation} 
Based on (\ref{aa44})--(\ref{aaa47})  and (\ref{aa37}) we conclude that 
\begin{equation}
\label{aa48}
Z_{[\xi_j,  \xi_j]}(\Psi_{\alpha}(x,  \xi_j)) =  Z_{I_k}(\Psi_{\alpha}(x, \xi_j))=1  
\ \ \mbox{for any}\ \  k=2,3,\ldots, n-1,
\end{equation} 
and 
\begin{equation}
\label{aaa48}
Z_{I_1}(\Psi_{\alpha}(x, \xi_j))=Z_{I_n}(\Psi_{\alpha}(x, \xi_j))=0.
\end{equation}
Thus, the first statement of Lemma~\ref{lemma5} is also valid for $\nu=\xi_j, \  
j=1, 2, \ldots, n-1.$

Now assume that $\alpha<-(n-1).$ Let $\nu \in I_j=(\xi_{j-1}, \xi_j), \ j=2,3,\ldots, n-1.$
Applying the first statement of Lemma~\ref{lemma4} in the case of finite intervals $I_k,$ 
and additionally  the third statement of Lemma~\ref{lemma4} in the case of infinite intervals 
$I_1$ or $I_n,$ taking into 
account (\ref{aa37}) we obtain
\begin{equation}
\label{ss1}
Z_{I_k}(\Psi_{\alpha}(x, \nu))=1,  \  \  k \neq j,
\end{equation} 
while
\begin{equation}
\label{ss2}
Z_{I_j}(\Psi_{\alpha}(x, \nu))=0.
\end{equation}
Using the same idea, in the case of $\nu \geq \xi_{n-1}$ we have
\begin{equation}
\label{ss3}
Z_{I_k}(\Psi_{\alpha}(x, \nu))=1,  \  \  k =1,2,\ldots, n-1,
\end{equation} 
while
\begin{equation}
\label{ss4}
Z_{I_n}(\Psi_{\alpha}(x, \nu))=0;
\end{equation}
and in the case of $\nu \leq \xi_1$
\begin{equation}
\label{ss5}
Z_{I_k}(\Psi_{\alpha}(x, \nu))=1,  \  \  k =2,3,\ldots, n,
\end{equation} 
while
\begin{equation}
\label{ss6}
Z_{I_1}(\Psi_{\alpha}(x, \nu))=0.
\end{equation}

In the case of $\nu=\xi_j, j=2, 3, \ldots, n-2,$ we use again the first 
and third statements  of Lemma~\ref{lemma4} and (\ref{aa37}). 
We have
\begin{equation}
\label{ss7}
Z_{[\xi_j,  \xi_j]}(\Psi_{\alpha}(x,  \xi_j)) =  Z_{I_k}(\Psi_{\alpha}(x, \xi_j))=1, 
 \  \  k=1, 2, \ldots, j-1, j+1, \dots, n,
\end{equation} 
while
\begin{equation}
\label{ss8}
Z_{I_j}(\Psi_{\alpha}(x, \nu))=Z_{I_{j+1}}(\Psi_{\alpha}(x, \nu))=0.
\end{equation}
Lemma~\ref{lemma5} is proved. $\square$

\section{Proof of Theorem~\ref{Hawaiitrue} and Beyond} \label{Sec9}

As we noticed in the Section~\ref{Sec5} (see (\ref{b3})) there is an intimate connection between 
the functions $\varphi'_{\lambda}$ and $Q_{\varkappa}[p].$ The following lemma 
gives fundamental estimations for the proof of Theorem~\ref{Intend}.
\begin{lemma}\label{lemma6} The following statements are true.
\begin{enumerate}
\item If $\lambda > -\frac{1}{n},$ then
\begin{itemize}
\item   For any interval $I_k$ of the second type
\begin{equation}
\label{tt1}
Z_{I_k}(\varphi_{\lambda}'(x))=1.
\end{equation}
\item  For any interval $I_k$ of the first type
\begin{equation}
\label{tt2}
Z_{I_k}(\varphi_{\lambda}'(x))=0.
\end{equation}
\end{itemize}
\item  If $\lambda\leq -1,$ then 
\begin{itemize}
\item For any finite interval of the second type $I_k$ 
\begin{equation}
\label{tt3}
 Z_{I_k}(\varphi_{\lambda}'(x))=1.
\end{equation}
\item For any infinite interval of the first type $I_k$ 
\begin{equation}
\label{tt4}
Z_{I_k}(\varphi_{\lambda}'(x))=1.
\end{equation}
\item For any infinite interval of the second type $I_k$ 
\begin{equation}
\label{tt5}
Z_{I_k}(\varphi_{\lambda}'(x))=0.
\end{equation}
\end{itemize}
\end{enumerate}
\end{lemma}

{\bf Proof of Lemma~\ref{lemma6}.}  By virtue of (\ref{aa24}) and (\ref{aa28}) the statement 
of Lemma~\ref{lemma5} can be rewritten as follows.
\begin{lemma}\label{lemma5'} The following statements are true.
\begin{enumerate}
\item If $\lambda > 0$ or $\lambda \leq -1,$ then for all real values of the parameter $\mu$ we have 
\begin{itemize}
\item $Z_{I_k}(\Phi_{\lambda}'(x, \mu))\leq 2$ when $k=2, 3, \ldots, n-1.$
\item $Z_{I_1}(\Phi_{\lambda}'(x, \mu))\leq 1$ and $Z_{I_n}(\Phi_{\lambda}'(x, \mu))\leq 1.$
\end{itemize}
\item  If $-\frac{1}{n}< \lambda <0,$ then for all real values of the parameter $\mu$ and all 
$k=1, 2, \ldots, n$ we have
 $$Z_{I_k}(\Phi_{\lambda}'(x, \mu))\leq 1.$$
\end{enumerate}
\end{lemma}

Let us apply Rolle's theorem to the above statements. If $\lambda > 0$ or $\lambda \leq -1,$ 
then for all real values of the parameter $\mu$ we have
\begin{equation}
\label{ss10}
 Z_{I_k}(\Phi_{\lambda}(x, \mu))\leq 3, \ \  k=2, 3, \ldots, n-1,
\end{equation}
 and
\begin{equation}
\label{ss11}
Z_{I_1}(\Phi_{\lambda}(x, \mu))\leq 2, \ \ Z_{I_n}(\Phi_{\lambda}(x, \mu))\leq 2.
\end{equation}
 Assume that $-\frac{1}{n}< \lambda <0.$ Then for all real values of the parameter $\mu$ we have
\begin{equation}
\label{ss12}
Z_{I_k}(\Phi_{\lambda}(x, \mu))\leq 2, \ \ k=1, 2, \ldots, n.
\end{equation}

Formulas (\ref{ss10}), (\ref{ss11}), (\ref{ss12})  together with Lemma~\ref{lemma1} and 
(\ref{aa9}) provide the following facts. 

\begin{enumerate}
\item{\it If $\lambda > 0,$ and  $\mu$ is any real number, then}
\begin{itemize} 
\item{\it For any finite interval $I_k$ of the second type}
\begin{equation}
\label{ss13}
 Z_{I_k}(\Phi_{\lambda}(x, \mu))= Z_{I_k}(\varphi_{\lambda}(x)- \mu)=0 \ \mbox{or} \  2.
\end{equation}
\item{\it For any infinite interval $I_k$  of the first type  }
\begin{equation}
\label{tt6}
Z_{I_k}(\Phi_{\lambda}(x, \mu))=Z_{I_k}(\varphi_{\lambda}(x)- \mu)=1.
\end{equation}
\item{\it For any infinite interval  $I_k$  of the second type }
\begin{equation}
\label{tt7}
Z_{I_k}(\Phi_{\lambda}(x, \mu))=Z_{I_k}(\varphi_{\lambda}(x)- \mu)=0\ \mbox{or} \ 2.
\end{equation}
\end{itemize}
\item{\it  If $ -\frac{1}{n}< \lambda < 0,$ and $\mu$ is any real number, then } 
\begin{itemize}
\item{\it For any interval $I_k, \ \ k=1,\ldots, n,$ of the first type}
\begin{equation}
\label{tt8}
Z_{I_k}(\Phi_{\lambda}(x, \mu))=Z_{I_k}(\varphi_{\lambda}(x)- \mu)=1.
\end{equation}
\item{\it For any interval $I_k, \ \ k=1, \ldots, n,$ of the second type}
\begin{equation}
\label{ttt8}
Z_{I_k}(\Phi_{\lambda}(x, \mu))=Z_{I_k}(\varphi_{\lambda}(x)- \mu)=0 \ \mbox{or} \  2.
\end{equation}
\end{itemize}
\item{\it If $  \lambda \leq -1,$ and $\mu$ is any real number, then  } 
\begin{itemize}
\item{\it For any finite interval $I_k$ of the second type}
\begin{equation}
\label{sss13}
 Z_{I_k}(\Phi_{\lambda}(x, \mu))= Z_{I_k}(\varphi_{\lambda}(x)- \mu)=0 \ \mbox{or} \  2.
\end{equation}
\item{\it For any infinite interval $I_k$ of the first type  }
\begin{equation}
\label{tt9}
Z_{I_k}(\Phi_{\lambda}(x, \mu))=Z_{I_k}(\varphi_{\lambda}(x)- \mu)=0 \ \mbox{or} \ 2.
\end{equation}
\item {\it  For any infinite interval $I_k$ of the second type  }
\begin{equation}
\label{ttt6}
Z_{I_k}(\Phi_{\lambda}(x, \mu))=Z_{I_k}(\varphi_{\lambda}(x)- \mu)=1.
\end{equation}
\end{itemize}
\end{enumerate}

Taking into account  (\ref{a2}), (\ref{a3}), (\ref{c1}), and applying the second statement of 
Lemma~\ref{phi0} in the case of $\lambda=0,$ we come to the conclusion.
\begin{enumerate}
\item{\it If $\lambda \geq 0,$}
\begin{itemize} 
\item{\it For any finite interval $I_k$ of the second type}
\begin{equation}
\label{sst13}
 Z_{I_k}(\varphi_{\lambda}'(x))=1.
\end{equation}
\item{\it For any infinite interval $I_k$  of the first type  }
\begin{equation}
\label{tts6}
Z_{I_k}(\varphi_{\lambda}'(x))=0.
\end{equation}
\item{\it For any infinite interval $I_k$  of the second type }
\begin{equation}
\label{tts7}
Z_{I_k}(\varphi_{\lambda}'(x))=1.
\end{equation}
\end{itemize}
\item{\it If $ -\frac{1}{n}< \lambda < 0,$ } 
\begin{itemize}
\item{\it For any interval $I_k, \ \ k=1, \ldots, n,$ of the first type}
\begin{equation}
\label{tts8}
Z_{I_k}(\varphi_{\lambda}'(x))=0.
\end{equation}
\item{\it For any interval $I_k, \ \ k=1, \ldots, n,$ of the second type}
\begin{equation}
\label{ttts8}
Z_{I_k}(\varphi_{\lambda}'(x))=1.
\end{equation}
\end{itemize}
\item{\it If $  \lambda \leq -1,$ }
\begin{itemize}
\item{\it For any finite interval $I_k$ of the second type}
\begin{equation}
\label{ssst13}
 Z_{I_k}(\varphi_{\lambda}'(x))=1.
\end{equation}
\item{\it For any infinite interval $I_k$  of the first type }
\begin{equation}
\label{tts9}
Z_{I_k}(\varphi_{\lambda}'(x))=1.
\end{equation}
\item {\it  For any infinite interval $I_k$  of the second type   }
\begin{equation}
\label{ttts6}
Z_{I_k}(\varphi_{\lambda}'(x))=0.
\end{equation}
\end{itemize}
\end{enumerate}

In order to complete the proof of Lemma~\ref{lemma6} we only have to check (\ref{tt2}) in 
the case of $\lambda >0.$ It follows from the first statement of Lemma~\ref{phi0} that if 
$\lambda>0,$ the function $$\varphi_{\lambda}(x))=\frac{p(x)}{p'(x)}+\lambda x $$
is increasing on any interval $I_k, \ k=1, 2, \ldots, n,$ of the first type. Hence, on any interval 
$I_k$ of the first type
\begin{equation}
\label{att1}
Z_{I_k}(\varphi_{\lambda}'(x))=0.
\end{equation}
Lemma~\ref{lemma6} is proved. $\square$

Let us prove Theorem~\ref{Hawaiitrue}. We have $\varkappa>\frac{n-1}{n}.$ According to (\ref{b3})
\begin{equation}\label{olga}
Z_{I_k}(\varphi_{\lambda}'(x))=Z_{I_k}\left(Q_{\varkappa }[p](x)\right),
\end{equation}
where by (\ref{b2}) 
$$\lambda=\varkappa-1>-\frac{1}{n}.$$
Now we can rewrite the first statement of Lemma~\ref{lemma6} as follows.
Let $\varkappa>\frac{n-1}{n}$. Then
 \begin{equation}
\label{II2}
Z_{I_k}\left( Q_{\varkappa}[p]\right)=1, \quad \mbox {\it if} \ I_k \ \mbox{\it is any interval of the second type}.
\end{equation}
 \begin{equation}
\label{II1}
Z_{I_k}\left( Q_{\varkappa}[p]\right)=0 \quad \mbox{\it if} \ I_k \ \mbox{\it is any interval of the first type}.
\end{equation}
Thus,  
$Z_{\mathbb{R}}\left(Q_{\varkappa}[p]\right)$ equals to the number of all intervals of the 
second type, that is
$$Z_{\mathbb{R}}\left(H_{\varkappa }[p](x)\right)=Z_{\mathbb{R}}
\left(Q_{\varkappa }[p](x)\right)=Z_{\mathbb{R}}(\varphi_{\lambda}'(x))=
n-Z_{\mathbb{R}}(p)=Z_{\mathbb{C}}(p).$$
Theorem~\ref{Hawaiitrue} is proved. $\square$

Lemma~\ref{lemma6} allows to make important conclusions in the case of $\varkappa\leq 0.$
If $\varkappa\leq 0,$ then $\lambda=\varkappa-1\leq -1.$ The following lemma follows from 
(\ref{ssst13}), (\ref{tts9}) and (\ref{ttts6}).

\begin{lemma}\label{lemmanew} If $\varkappa\leq 0,$ then
\begin{itemize}
\item For any finite interval $I_k$ of the second type
\begin{equation}
\label{sstt13}
 Z_{I_k}\left(H_{\varkappa }[p](x)\right)=Z_{I_k}\left(Q_{\varkappa }[p](x)\right)=1.
\end{equation}
\item For any infinite interval $I_k$  of the first type
\begin{equation}
\label{tsss90}
Z_{I_k}\left(H_{\varkappa }[p](x)\right)=Z_{I_k}\left(Q_{\varkappa }[p](x)\right)=1.
\end{equation}
\item For any infinite interval $I_k$  of the second type 
\begin{equation}
\label{ttss6}
Z_{I_k}\left(H_{\varkappa }[p](x)\right)=Z_{I_k}\left(Q_{\varkappa }[p](x)\right)=0.
\end{equation}
\end{itemize}
\end{lemma}

\begin{remark} \label{kappal0}
Note that in order to prove Theorem~\ref{Theorem2} we just have to show that for 
any finite interval $I_k$ of the first type 
\begin{equation}
\label{ttss7}
Z_{I_k}\left( H_{\varkappa}[p]\right)=Z_{I_k}(\varkappa (p')^2-p''p)=2.
\end{equation}
\end{remark}
Indeed, if we use the standard notation $\#S$ for the number of elements in any finite set $S,$  we obtain
$$Z_{\mathbb{R}}\left(H_{\varkappa }[p](x)\right)=Z_{\mathbb{R}}\left(Q_{\varkappa }[p](x)\right)=$$  
$$\#\{I_k | \ I_k \ \mbox{is a finite interval of the second type}\}+$$
$$+2\cdot \# \{I_k| \ I_k \ \mbox{is a finite interval of the first type}\}
+ $$  $$\#\{I_k| \ I_k \ \mbox{is an infinite interval of the first type}\}$$
$$= \# \{I_k| \ I_k \ \mbox{is a finite interval}\}
+\#\{I_k| \ I_k \ \mbox{is an interval of the first type}\}=$$  $$n-2+Z_{\mathbb{R}}(p).$$

\section{ Relations Between Zeros of $H_{\varkappa }[p](x)$ and Zeros of $H_{2-\frac{1}{\varkappa}}[p'](x).$ 
Proof of Theorem~\ref{Theorem2}} \label{Sec10}

Using definition (\ref{main.function.2.numerator}) of the function $H_{\varkappa }[p]$ we 
obtain the following formula for the derivative $H_{\varkappa }'[p]$
\begin{equation}
\label{ttss9}
H_{\varkappa}'[p]=(2\varkappa-1) p' p''-p'''p,
\end{equation}
and if $\varkappa \neq 0,$ then
\begin{equation}
\label{ttss10}
H_{2-\frac{1}{\varkappa}}[p']=\left(2-\frac{1}{\varkappa}\right) (p'')^2-p'''p'. 
\end{equation}

The lemma below can be easily verified.
\begin{lemma}\label{lemma7} For $\varkappa \neq 0,$
\begin{equation}
\label{ttss11}
H_{\varkappa}'[p]\cdot p'-H_{2-\frac{1}{\varkappa}}[p']\cdot p=
\frac{2\varkappa-1}{\varkappa} H_{\varkappa}[p]\cdot p'',
\end{equation}
\begin{equation}
\label{ttss12}
H_{\varkappa}'[p]\cdot p''-\varkappa H_{2-\frac{1}{\varkappa}}[p']\cdot p'= H_{\varkappa}[p]\cdot p'''.
\end{equation}
\end{lemma}
The statement below can be considered as an analogue of Rolle's Theorem, which provides the relation between zeros of $H_{\varkappa}[p]$ and zeros of $H_{\left(2-\frac{1}{\varkappa}\right)}[p'].$
\begin{lemma}\label{lemma8} Assume that $\varkappa \neq 0,$ and
\begin{enumerate}
\item{\it $H_{\varkappa}[p](x) \neq 0$ for all $x \in (a,b).$}
\item $H_{\varkappa}[p](a) =H_{\varkappa}[p](b)=0.$
\end{enumerate}
Then
\begin{enumerate}
\item{\it Let $p'(x) \neq 0$ and $p(x)\neq 0$ for all $x \in [a,b].$ There 
exists $y \in (a,b)$ such that $$H_{2-\frac{1}{\varkappa}}[p'](y)=0.$$}
\item {\it Let $p'(x) \neq 0$ and $p''(x)\neq 0$ for all $x \in [a,b].$ There exists 
$y \in (a,b)$ such that $$H_{2-\frac{1}{\varkappa}}[p'](y)=0.$$}
\end{enumerate}
\end{lemma}

{\bf Proof of Lemma~\ref{lemma8}} Let us prove the first statement of Lemma~\ref{lemma8}. 
Assume that $p(x)\neq 0$ and $p'(x)\neq 0$ for all $x \in (a,b).$ Let $a$ be a root of multiplicity 
$m$ and $b$ be a root of multiplicity $k$ of the function $H_{\varkappa}[p].$ It means that there 
are polynomials $q_1$ and $q_2$ as well as constants $c\neq 0$ and $d\neq 0$ such that
$$H_{\varkappa}[p](x)=(x-a)^m(c+(x-a)q_1(x)), \quad \mbox{and}$$
\begin{equation}\label{EQ1}
   H_{\varkappa}[p](x)=(x-b)^k(d+(x-b)q_2(x)).
\end{equation}
Because $H_{\varkappa}[p](x)\neq 0$ in $(a,b),$ we have
\begin{equation}\label{EQ2}
(-1)^k c\cdot d>0.
\end{equation}
Since $p(x)\neq 0$ and $p'(x)\neq 0$ it follows from (\ref{ttss11}), (\ref{ttss12}) that $a$ 
is a root of multiplicity at least $m-1$ and $b$ is a root of multiplicity at least $k-1$ of the function 
$H_{\left(2-\frac{1}{\varkappa}\right)}[p'].$ So, there are polynomials $T_1$ and $T_2$ such that
\begin{equation}\label{EQ3}
H_{2-\frac{1}{\varkappa}}[p'](x)=(x-a)^{m-1}T_1(x) \quad \mbox{and} \quad 
H_{2-\frac{1}{\varkappa}}[p'](x)=(x-b)^{k-1}T_2(x).
\end{equation}
Using (\ref{ttss11}), (\ref{EQ1}) and (\ref{EQ3}) we obtain
\begin{equation}
\label{ttss13}
 T_1(a)=\frac{p'(a)}{p(a)} m c,
\end{equation}
\begin{equation}
\label{ttss14}
 T_2(b)=\frac{p'(b)}{p(b)} k d.
\end{equation}
Because $p'(x)\neq 0$ and $p(x) \neq 0$ on the interval $[a,b],$ the following is true
\begin{equation}
\label{ttss15}
\frac{p'(a)}{p(a)} \cdot \frac{p'(b)}{p(b)}>0.
\end{equation}
Therefore, by (\ref{EQ2})
$$(-1)^k T_1(a)T_2(b)>0.$$
Since the functions $T_1(x)$ and $T_2(x)$ are continuous, there exists $\varepsilon>0$ such that
\begin{equation}\label{ttss16}
(-1)^kT_1(a+\varepsilon)T_2(b-\varepsilon)>0,
\end{equation}
and by (\ref{EQ3}) 
\begin{equation}
\label{ttss17}
H_{2-\frac{1}{\varkappa}}[p'](a+\varepsilon)\cdot H_{2-\frac{1}{\varkappa}}
[p'](b-\varepsilon)=(\varepsilon)^{m+k-2}(-1)^{k-1}T_1(a+\varepsilon)T_2(b-\varepsilon)<0.
\end{equation}
We obtain the first statement of Lemma~\ref{lemma8}.

The second statement of Lemma~\ref{lemma8} can be derived from (\ref{ttss12}) in the similar way.
Lemma~\ref{lemma8} is proved. $\square$

The following estimation for the number $Z_{(a,b)}\left(H_{\varkappa}[p] \right)$ through 
$Z_{(a,b)}\left(H_{2-\frac{1}{\varkappa}}[p'] \right)$ is  derived from Lemma~\ref{lemma8} 
in the same way as similar facts about relations between the number of zeros of a function $f$ 
and the number of zeros of its derivative $f'$ are usually obtained from the original Rolle's theorem. 
\begin{corol}
\label{Rolle}
Assume that $\varkappa \neq 0.$  Let either $p'(x)\neq 0$ and $p(x)\neq 0$  
for all $x\in [a,b],$ or $p''(x)\neq 0$ and $p'(x)\neq 0$  for all $x\in [a,b].$
Then
\begin{equation}\label{R1}
Z_{(a,b)}\left(H_{\varkappa}[p] \right)\ \leq \ Z_{(a,b)}
\left(H_{2-\frac{1}{\varkappa}}[p'] \right)+1.
\end{equation}
\end{corol}
Denote by $\gamma_1 <\gamma_2< \ldots <\gamma_{n-2}$ the roots of $p''.$ Since the 
roots $\xi_1 <\xi_2< \ldots <\xi_{n-1}$ of the polynomial $p'$ are real and simple, we have
\begin{equation}
\label{ttss21}
-\infty<\xi_1 <\gamma_1 <\xi_2<\gamma_2< \xi_3<\ldots <\gamma_{n-2}<\xi_{n-1}<\infty.
 \end{equation}
We will consider the intervals
$$J_{2k-1}=(\xi_k,\gamma_k), \ J_{2k}=( \gamma_k, \xi_{k+1}),\ \ k=1,2,\ldots, n-2;$$
\begin{equation}
\label{ttss19}
J_0=I_1=(-\infty,\xi_1), \  \  J_{2n-3}=(\xi_{n-1}, \infty).
\end{equation} 
\begin{remark}\label{intJ}
Both $p'$ and $p''$ do not vanish on each of the intervals $J_k,  \  k=0,1, \ldots, 2n-3.$
\end{remark}
The following fact about the number of zeros of the function $H_{\varkappa}[p]$ on 
the intervals $J_k$ is true.
\begin{lemma}\label{lemma9} Let $\varkappa<0.$ Then 
\begin{equation}
\label{ttss22}
Z_{J_m}(H_{\varkappa}[p])\leq 1
\end{equation}
for any $m=0,1,\ldots, 2n-3.$
\end{lemma}

{\bf Proof of Lemma~\ref{lemma9}} Given $\varkappa<0,$ we have $2-\frac{1}{\varkappa}>2.$ 
Since all roots of $p'$ are real and simple, all roots of its derivative $p''$ are also real and simple. 
Therefore, we can apply Theorem~\ref{Hawaiitrue} to the polynomial $p'.$ We have
\begin{equation}
\label{ttss18}
Z_{\mathbb{R}}\left(H_{2-\frac{1}{\varkappa}}[p']\right)=Z_{\mathbb{C}}(p')=0.
\end{equation}

Let us fix any interval $J_m, \ m=0,1,\ldots, 2n-3.$
Because by Remark~\ref{intJ} both $p'$ and $p''$ do not vanish on $J_m,$ 
we can apply Corollary~\ref{Rolle} and obtain the statement of Lemma\ref{lemma9}.
Lemma\ref{lemma9} is proved. $\square$

The following lemma is the key to the proof of Theorem~\ref{Theorem2}.
\begin{lemma}\label{lemma10}  Let $\varkappa\leq 0.$ Then
\begin{equation}
\label{ttss26}
Z_{I_k}(H_{\varkappa}[p])\leq 2
\end{equation}
for any $k=2,3,\ldots, n-1,$ and
\begin{equation}
\label{ttss262}
Z_{I_1}(H_{\varkappa}[p])\leq 1, \ \ Z_{I_n}(H_{\varkappa}[p])\leq 1.
\end{equation}
\end{lemma}

{\bf Proof of Lemma~\ref{lemma10}}.
By (\ref{ttss21}) and (\ref{ttss19})
$$I_k=(\xi_{k-1}, \xi_k)=(\xi_{k-1}, \gamma_{k-1})\cup \{\gamma_{k-1}\} 
\cup ( \gamma_{k-1}, \xi_k)=J_{2k-3}\cup \{\gamma_{k-1}\} \cup J_{2k-2}.$$  
Note that
\begin{equation}
\label{ttss25}
H_{\varkappa}[p](\gamma_k)\neq 0, \  \ k=1,2,\ldots, n-2.
\end{equation}
Indeed, suppose that $H_{\varkappa}[p](\gamma_k)= 0.$ Since $\gamma_k$ is a root of 
$p''$, it follows from the definition (\ref{main.function.2.numerator}) of the function 
$H_{\varkappa}[p]$ that in this case $p'(\gamma_k)=0$. This fact contradicts to the assumption 
of Theorem~\ref{Theorem2} that $p'$ doesn't have multiple roots. 

If $\varkappa<0$ the statement of lemma follows immediately from (\ref{ttss22}) and  (\ref{ttss25}). 
Note that for $\varkappa=0$ we have
$$H_{\varkappa}[p](x)=-p(x)p''(x).$$
Since all roots of $p'$ are real and simple, there is exactly one root of $p''$ inside each finite interval $I_k.$ 
If $I_k$ is an interval of the first type, it also contains one root of $p,$  while if it is of the second type it 
doesn't have roots of $p.$ So,
\begin{equation}\label{I10finite}
Z_{I_k}(H_{\varkappa}[p])= 2,\quad \mbox{\it if}  \ I_k \ \mbox{\it is a finite interval of the first type,}
\end{equation}
and
\begin{equation}\label{I20finite}
Z_{I_k}(H_{\varkappa}[p])= 1, \quad \mbox{\it if}  \ I_k \ \mbox{\it is a finite interval of the second type.}
\end{equation}
If $I_k$ is an infinite interval, there is no roots of $p''$ in $I_k.$ Thus,
\begin{equation}\label{I10infinite}
Z_{I_k}(H_{\varkappa}[p])= 1, \ \ \mbox{\it if}  \ I_k  \ \mbox{\it is an infinite interval of the first type,}
\end{equation}
and
\begin{equation}\label{I20infinite}
Z_{I_k}(H_{\varkappa}[p])= 0, \ \ \mbox{\it if}  \ I_k  \ \mbox{\it is an infinite interval of the second type.}
\end{equation}
Lemma\ref{lemma10} is proved.  $\square$ 

Since $p$ doesn't have real multiple roots, and all roots of $p'$ are real and simple,  
\begin{equation}
\label{ttss50}
H_{\varkappa}[p](x)=0 \ \Leftrightarrow  \ Q_{\varkappa}[p](x)=0.
\end{equation}
So, for each $k=1,2,\ldots, n$ and any real $\varkappa$ we have 
\begin{equation}
\label{ttss502}
Z_{I_k} \left(H_{\varkappa}[p]\right)=Z_{I_k} \left(Q_{\varkappa}[p]\right).
\end{equation}
Now for each $\varkappa\leq 0$ formula (\ref{ttss7}) follows from Lemma~\ref{lemma10} 
and Lemma~\ref{IntI1}, that is
\begin{equation}\label{I1finite}
Z_{I_k}(H_{\varkappa}[p])= 2,\quad \mbox{\it if}  \ I_k \ \mbox{\it is a finite interval of the first type}.
\end{equation}
By virtue of Remark~\ref{kappal0} this completes the proof of Theorem~\ref{Theorem2}.

\section{ Some properties of $M[p]$} \label{Sec11}

 Assume that $I_k$ is a finite interval of the first type. According to (\ref{I1finite}) for each 
 $\varkappa\leq 0$ there are exactly two roots of $Q_{\varkappa}[p](x)$ inside $I_k.$ 
 The statement below gives an idea about the zero distribution of the function 
 $Q_{\varkappa}[p](x)$ inside a finite interval $I_k$ of the first type for $\varkappa<0.$
\begin{lemma}\label{lemma11}  Let $\varkappa<0.$ Suppose that $I_k=(\xi_{k-1}, \xi_k)$ is 
a finite interval of the first type, and $\alpha \in I_k$ is the only zero of $p$ inside the interval 
$I_k.$ Then
\begin{equation}
\label{sstt39}
Z_{(\xi_{k-1}, \alpha)} \left(Q_{\varkappa}[p]\right)=Z_{(\alpha,\xi_k)} \left(Q_{\varkappa}[p]\right)=1.
\end{equation}
\end{lemma}

{\bf Proof of Lemma~\ref{lemma11}} Note that by the definition (\ref{ttss27}) of the function $M[p]$ we have
$$M[p](\alpha)=0,$$
and by (\ref{ttss28}) for any $\varkappa<0$ we have
\begin{equation}
\label{sstt40}
Q_{\varkappa}[p](\alpha)=\varkappa \ -M[p](\alpha)<0.
\end{equation}
Using (\ref{IRR}) and (\ref{ILR}) we obtain
\begin{equation}
\label{IRRQ}
\lim_{x\to\xi_{k-1}^{+}}Q_{\varkappa}[p](x)=\varkappa \ -\lim_{x\to\xi_{k-1}^{+}}M[p])(x)>0,
\end{equation}
and
\begin{equation}
\label{ILRQ}
\lim_{x\to\xi_{k}^{-}}Q_{\varkappa}[p](x)= \varkappa-\lim_{x\to\xi_{k}^{-}} M[p](x)>0.
\end{equation}
It follows from (\ref{IRRQ}), (\ref{ILRQ}) and (\ref{sstt40}) that both numbers $Z_{(\xi_{k-1}, \alpha)} 
\left(Q_{\varkappa}[p]\right)$ and $Z_{(\alpha,\xi_k)} \left(Q_{\varkappa}[p]\right)$ are odd. By virtue 
of (\ref{I1finite}) it is possible only if the statement of Lemma~\ref{lemma11} is true. $\square$

\begin{remark}\label{gamma}
Let $I_k=(\xi_{k-1}, \xi_k)$ be a finite interval of the first type, and $\gamma_{k-1} \in I_k$ be the 
only zero of $p''$ inside the interval $I_k.$ In the same way as above, using (\ref{IRRQ}), (\ref{ILRQ}), 
(\ref{I1finite}), (\ref{ttss22}), and (\ref{ttss25}) we can prove that
 \begin{equation}
\label{SSTT39}
Z_{(\xi_{k-1}, \gamma)} \left(Q_{\varkappa}[p]\right)=Z_{(\gamma, \ \xi_k)} 
\left(Q_{\varkappa}[p]\right)=1.
\end{equation}
\end{remark}
\begin{remark}\label{graph}
The relations (\ref{sstt39}) and (\ref{SSTT39}) give an idea how the graph of the function 
$M[p](x)$ should look on finite intervals of the first type. It follows from (\ref{ttss27}) that the 
function $M[p]$ has exactly two roots in the interval $I_k$: a root $\alpha$ of $p$ and  a root 
$\gamma$ of $p''.$  Let us denote these roots of $M[p]$ by $z_1$ and $z_2,$ where $z_1 \leq z_2.$ 
By virtue of (\ref{ttss28}) the relations (\ref{sstt39}) and (\ref{SSTT39}) mean that for each 
$\varkappa \leq 0$ the equation (\ref{mainequation})
$$ M[p](x)=\varkappa$$ 
has the unique solution in the interval $(\xi_{k-1}, z_1)$ and the unique solution in the interval 
$(z_2, \xi_k).$ Taking into account (\ref{IRR}) and (\ref{ILR}) we can conclude that
\begin{equation}\label{GrM1}
M[p](x)< 0 \quad \mbox{\it and is increasing on the interval} \ (\xi_{k-1}, z_1).
\end{equation}
\begin{equation}\label{GrM2}
M[p](x)< 0 \quad \mbox{\it and is decreasing on the interval} \ (z_2, \xi_k).
\end{equation}
If additionally we use (\ref{II1}), we can state that 
$$M[p](x)>0, \quad x \in (z_1,z_2),$$
and
$$\max_{x\in [z_1, z_2]} M[p](x)<\frac{n-1}{n}.$$
\end{remark}

\section{ The case of $0< \varkappa <\frac{1}{2} $} \label{Sec12}
In this section we will prove the following statement. Theorem~\ref{Theorem3} is a direct corollary 
from this result.
\begin{theorem}\label{Theorem4}
 Let $0\leq \varkappa<\frac{1}{2}.$ Then 
\begin{enumerate}
\item  $Z_{I_k}\left(Q_{\varkappa}[p]\right)=Z_{I_k}\left(H_{\varkappa}[p]\right)=1,$ if 
$I_k$ is a finite interval of the second type.
\item $Z_{I_k}\left(Q_{\varkappa}[p]\right)=Z_{I_k}\left(H_{\varkappa}[p]\right) =0 \ \mbox{or} 
\ 2,$ if $I_k$ is a finite interval of the first type.
\item $Z_{I_k}\left(Q_{\varkappa}[p]\right)=Z_{I_k}\left(H_{\varkappa}[p]\right) =0,$  if $I_k$ 
is an infinite interval of the second type.
\item $Z_{I_k}\left(Q_{\varkappa}[p]\right)=Z_{I_k}\left(H_{\varkappa}[p]\right) =1,$  if $I_k$ 
is an infinite interval of the first type.
\end{enumerate}
\end{theorem}

{\bf Proof of Theorem~\ref{Theorem4} }
If $\varkappa=0,$ then $Q_{\varkappa}[p]=M[p]$ and by definition (\ref{ttss27}) of the function 
$M[p]$ all statements of theorem are obvious.

Now assume that $\varkappa\neq 0$.
Let $I_k$ be a finite interval of the first type.
Note that one of the numbers $z_1$ and $z_2,$ described in Remark~\ref{GrM2}, is a root 
$\gamma_{k-1}$ of the polynomial $p''.$ Hence,
\begin{equation}\label{EQU4}
  (z_1,z_2)\subset (\xi_{k-1}, \gamma_{k-1})=J_{2k-3} \quad \mbox{or} \quad  
  (z_1,z_2)\subset ( \gamma_{k-1}, \xi_k)=J_{2k-2}. 
\end{equation}
It follows from Remark~\ref{GrM2} that the equation (\ref{mainequation}) doesn't have solutions in 
the intervals $(\xi_{k-1}, z_1)$ and $(z_2, \xi_k).$  Let us choose the interval from $J_{2k-3}$ and 
$J_{2k-2}$ that contains $(z_1,z_2),$ and denote it by $J.$ We have
\begin{equation}\label{EQU5}
Z_{I_k}\left(Q_{\varkappa}[p]\right)=Z_{(z_1,z_2)}\left(Q_{\varkappa}[p]\right)=
Z_J\left(Q_{\varkappa}[p]\right),
\end{equation}
if $I_k$ is a finite interval of the first type.

If $I_k=(\xi_{k-1},\xi_k)$ is a finite interval of the second type, it contains only one root of the 
function $M[p]$: the root $\gamma_{k-1}$ of the polynomial $p''.$ Let us choose the interval from 
$(\xi_{k-1}, \gamma_{k-1})$ and $(\gamma_{k-1}, \xi_k),$ where $M[p](x)>0.$ Now denote this 
interval by $J.$ If $0<\varkappa<\frac{1}{2},$ and the equation
$M[p](x)=\varkappa$
has a solution in $I_k,$ it must belong to $J.$ So,
\begin{equation}\label{EQU6}
Z_{I_k}\left(Q_{\varkappa}[p]\right)=Z_J\left(Q_{\varkappa}[p]\right),
\end{equation}
if $I_k$ is a finite interval of the second type.

If $0<\varkappa<\frac{1}{2},$ then $-\infty<2-\frac{1}{\varkappa}<0.$ Since all roots 
$(\xi_k)_{k=1}^{n-1}$ of $p'$ are real and simple, all roots $(\gamma_k)_{k=1}^{n-2}$ 
of $p''$ are also real and simple. So, $p'$ satisfies the hypotheses of Theorem~\ref{Theorem2}. 
The intervals
 $$(-\infty, \gamma_1), (\gamma_1, \gamma_2), \ldots, (\gamma_{n-3}, \gamma_{n-2}), 
 (\gamma_{n-2}, \infty)$$ 
are intervals (\ref{i1}) for $p'.$ All of them are of the first type. Because we denote by $\xi$ with indices 
the roots of $p'$ and $\gamma_{k-1}$ is  the root of $p'',$ the interval $J$ in both (\ref{EQU5}) and 
(\ref{EQU6}) plays exactly the same role for $p'$ as the intervals $(\xi_{k-1}, \alpha)$ and $(\alpha, \xi_k)$ 
for $p$ in Lemma~\ref{lemma11}. Applying Lemma~\ref{lemma11} to $p'$ on the interval $J$ we obtain
\begin{equation}\label{EQU7}
  Z_J\left(Q_{2-\frac{1}{\varkappa}}[p']\right)= Z_J\left(H_{2-\frac{1}{\varkappa}}[p']\right)  = 1.
\end{equation}
Note that $p'(x)\neq 0$ and $p''(x) \neq 0$ for all $x \in J.$ So, we can apply Corollary \ref{Rolle} 
from Lemma~\ref{lemma8} and obtain
\begin{equation}\label{EQU8}
  Z_J\left(H_{\varkappa}[p]\right)= Z_J\left(Q_{\varkappa}[p]\right)  \leq 2.
\end{equation}
Therefore, by (\ref{EQU5}) and  (\ref{EQU6}) for any finite interval $I_k$ we have
\begin{equation}\label{EQU9}
  Z_{I_k}\left(H_{\varkappa}[p]\right)= Z_{I_k}\left(Q_{\varkappa}[p]\right)  \leq 2.
\end{equation}
According to Lemma~\ref{Interval.endpoints} the number $Z_{I_k}\left(Q_{\varkappa}[p]\right) $ is even 
if $I_k$ is a finite interval of the first type, and it is odd if $I_k$ is a finite interval of the second type. We obtain 
first two statements of Theorem~\ref{Theorem4}.

Now consider the infinite intervals $I_1=(-\infty, \xi_1)$ and $I_n=(\xi_{n-1}, \infty).$ First, assume that an 
infinite interval $I$ is an interval of the first type. Because all root of $p'$ are real and simple, 
$$I_1\subset (-\infty,\gamma_1), \quad \mbox{and} \quad I_{n-1}\subset (\gamma_{n-2}, \infty).$$
Since $-\infty<2-\frac{1}{\varkappa}<0$, using the above remark by (\ref{ttss262}) we obtain
$$Z_{I}\left(H_{2-\frac{1}{\varkappa}}[p']\right)\leq 1. $$
The infinite intervals $I$ contain neither roots of $p'$ nor roots of $p''.$ We apply again 
Corollary~\ref{Rolle} from Lemma~\ref{lemma8} and obtain
$$Z_{I}\left(H_{\varkappa}[p]\right)= Z_{I}\left(Q_{\varkappa}[p]\right)  \leq 2.$$
Since $\varkappa<\frac{1}{2},$ by virtue of Lemma~\ref{Infint} if $I$ is an infinite interval of the 
first type, the number $Z_{I}\left(Q_{\varkappa}[p]\right) $ is odd. We conclude that
\begin{equation}\label{EQU10}
Z_{I}\left(H_{\varkappa}[p]\right)= Z_{I}\left(Q_{\varkappa}[p]\right)  =1.
\end{equation}
Let $I$ be an infinite interval of the second type. According to (\ref{ttss27}) 
\begin{equation}\label{EQU11}
M[p'](x)=\frac{p'''(x)p'(x)}{(p''(x))^2}.
\end{equation}
Because all roots of $p'$ are real and simple, it is clear that  $\xi_{n-1}$ is the only root of 
$M[p']$ on  $[\xi_{n-1}, +\infty)$, and $\xi_1$ is its only root on the interval $(-\infty,\xi_1].$ Since
$$\lim_{x\to \pm \infty} M[p'](x)=\frac{n-2}{n-1}>0,$$
we see that
\begin{equation}
\label{aabb2}
M[p'](x)>0, \quad \mbox{when} \ x\in (-\infty,\xi_1) \cup (\xi_{n-1}, +\infty).
\end{equation}
Since  $2-\frac{1}{\varkappa}<0,$  the equation
$$M[p'](x)=2-\frac{1}{\varkappa}$$
doesn't have solutions on infinite intervals of the second type. Thus, if $I$ is an infinite 
interval of the second type, then
\begin{equation}
\label{sttt42}
Z_{I} \left(Q_{2-\frac{1}{\varkappa}}[p']\right)=Z_{I} \left(H_{2-\frac{1}{\varkappa}}(p')\right)=0.
\end{equation}
Let us apply again Corollary~\ref{Rolle} from Lemma~\ref{lemma8}. We have
\begin{equation}
\label{sttt41}
Z_{I} \left(Q_{\varkappa}[p']\right)=Z_{I} \left(H_{\varkappa}(p')\right) \leq 1.
\end{equation}
According to Lemma~\ref{Infint} if $I$ is an infinite interval of the second type, the number 
$Z_{I}\left(Q_{\varkappa}[p]\right) $ is even. Therefore,
\begin{equation}\label{EQU12}
Z_{I}\left(H_{\varkappa}[p]\right)= Z_{I}\left(Q_{\varkappa}[p]\right)  =0.
\end{equation}
Theorem~\ref{Theorem4} is proved. $\square$

\section{ Some remarks about the intervals $I_k$ of the first type} \label{Sec13}

Because all roots of $p'$ are real and simple, the case of intervals of the first type is of the great 
interest for us. The following more general result can be obtained  in the same way as the second 
statement of Theorem~\ref{Theorem4}.
\begin{theorem}\label{Theorem5}
Let $0\leq \varkappa<\frac{n-1}{n}.$ Then for any finite interval of the first type $I_k$ the 
following is true
\begin{equation}
\label{abc1}
Z_{I_k}\left(Q_{\varkappa}[p]\right)=Z_{I_k}\left(H_{\varkappa}[p]\right) =0 \ \mbox{or} \ 2.
\end{equation}
\end{theorem}
This statement is equivalent to the following fact.
\begin{lemma}\label{lemma14}
 For every $j=1,2,\ldots,n-1,$ and $$\frac{j-1}{j}\leq \varkappa<\frac{j}{j+1},$$ on any 
finite interval $I_k$ of the first type the following estimation holds
\begin{equation}
\label{abc2}
Z_{I_k}\left(Q_{\varkappa}[p]\right)=Z_{I_k}\left(H_{\varkappa}[p]\right) =0 \ \mbox{or} \ 2.
\end{equation}
\end{lemma}

{\bf Proof of Lemma~\ref{lemma14}.} In order to prove Lemma~\ref{lemma14} we use induction 
on $j=1, 2, \ldots, n-1.$

The second statement of Theorem~\ref{Theorem4} provides the base of induction for $j=1.$
$$\mbox{If}\quad \frac{j-1}{j}\leq \varkappa<\frac{j}{j+1},\quad \mbox{then}\quad \frac{j-2}{j-1}
\leq 2-\frac{1}{\varkappa}<\frac{j-1}{j}. $$
Since all roots $(\xi_k)_{k=1}^{n-1}$ of $p'$ are real and simple, all roots $(\gamma_k)_{k=1}^{n-2}$ 
of $p''$ are also real and simple, and all the intervals $(\gamma_{k-1}, \gamma_k), \  k=2,3\ldots,n-2,$ 
are intervals of the first type for the polynomial $p'.$ Thus, by the induction hypothesis 
\begin{equation}
\label{abc3}
Z_{(\gamma_{k-1}, \gamma_k)}\left(H_{2-\frac{1}{\varkappa}}[p']\right)=0 \ \mbox{or} \ 2.
\end{equation}

Let us fix an interval $I_k=(\xi_{k-1}, \xi_k)$ of the first type. Denote by $\alpha$ the only root of the 
polynomial $p$ inside $I_k.$ Note that $\gamma_{k-1}$ is the only root of $p''$ inside $I_k.$ Without 
loss of generality we assume that
$$\xi_{k-1}<\gamma_{k-1} <\alpha <\xi_k$$
(the case when $\xi_{k-1}< \alpha<\gamma_{k-1} <\xi_k$ can be considered in the same way).
Obviously, we have the following inclusion
\begin{equation}\label{abc4}
(\gamma_{k-1}, \alpha) \ \subset \  (\gamma_{k-1}, \xi_k) \ \subset (\gamma_{k-1}, \gamma_k).
\end{equation}
It follows from(\ref{abc3}) that
\begin{equation}
\label{abc5}
Z_{(\gamma_{k-1}, \xi_k)}\left(H_{2-\frac{1}{\varkappa}}[p']\right)\leq 2.
\end{equation}
Let us apply Corollary \ref{Rolle} from Lemma~\ref{lemma8}. We obtain
\begin{equation}
\label{abc6}
Z_{(\gamma_{k-1}, \xi_k)}\left(H_{\varkappa}[p]\right)\leq 3.
\end{equation}
As we showed in Remark~\ref{GrM2}
\begin{equation}
\label{abc7}
Z_{(\xi_{k-1}, \xi_k)}\left(H_{\varkappa}[p]\right)=Z_{(\gamma_{k-1}, \alpha)}
\left(H_{\varkappa}[p]\right)=Z_{(\gamma_{k-1}, \xi_k)}\left(H_{\varkappa}[p]\right).
\end{equation}
Therefore,  
\begin{equation}
\label{abc8}
Z_{I_k}\left(Q_{\varkappa}[p]\right)=Z_{I_k}\left(H_{\varkappa}[p]\right)\leq 3.
\end{equation}
According to Lemma~\ref{Interval.endpoints}, the number $Z_{I_k}\left(Q_{\varkappa}[p]\right) $ 
is even if $I_k$ is a finite interval of the first type. Thus, for $\frac{j-1}{j}\leq \varkappa<\frac{j}{j+1}$ we obtain
\begin{equation}
\label{abc9}
Z_{I_k}\left(Q_{\varkappa}[p]\right)=Z_{I_k}\left(H_{\varkappa}[p]\right)\ = \ 0 \ \mbox{or} \ 2.
\end{equation}
We obtain the conclusion of Lemma~\ref{lemma14}. $\square$

\begin{remark}\label{refine}
Lemma~\ref{lemma14} allows to refine Remark~\ref{GrM2} about the behavior of the function 
$M[p]$ on a finite interval of the first type $I_k=(\xi_{k-1}, \xi_k)$: {\it there exists $x_{max} 
\in I_k$ such that}
\begin{equation}
\label{GrM3}
M[p](x) \quad \mbox{\it  is increasing on} \ (\xi_{k-1}, x_{max}],
\end{equation}
\begin{equation}
\label{GrM4}
M[p](x) \quad \mbox{\it is decreasing on} \ [x_{max}, \xi_k),
\end{equation}
{\it and}
\begin{equation}
\label{GrM5}
M[p]\left(x_{max}\right)<\frac{n-1}{n}.
\end{equation}
\end{remark}

The following statement gives the estimation for the number of real zeros of the function $H_{\varkappa}[p]$
on infinite intervals $I_k$ of the first type.
\begin{theorem}\label{Theorem6} If $0 \leq \varkappa<\frac{n-1}{n}, $ then for any infinite interval 
$I_k$ of the first type
\begin{equation}
\label{bbaa28}
Z_{I_k}(H_{\varkappa}[p])=Z_{I_k}\left(Q_{\varkappa}[p]\right) = 1.
\end{equation}
\end{theorem}
Theorem~\ref{Theorem6} follows from the lemma below. 

\begin{lemma}\label{lemma15}  For each $j=1,2,\ldots, n-1,$ and for every $$\frac{j-1}{j}
\leq \varkappa<\frac{j}{j+1}, \ \ j=1, 2, \ldots, n-1,$$ on any infinite interval of the first type  
$I_k$ the following estimation is true
\begin{equation}
\label{bbaa30}
Z_{I_k}\left(H_{\varkappa}[p]\right)=Z_{I_k}\left(Q_{\varkappa}[p]\right) = 1.
\end{equation}
\end{lemma}

{\bf Proof of Lemma~\ref{lemma15}.} We will use induction on $j.$
The fourth statement of Theorem~\ref{Theorem4} gives the base of induction for $j=1.$

$$\mbox{If}\quad \frac{j-1}{j}\leq \varkappa<\frac{j}{j+1},\quad \mbox{then}\quad \frac{j-2}{j-1}
\leq 2-\frac{1}{\varkappa}<\frac{j-1}{j}. $$
Given that $p'$ has only real and simple zeros, the induction hypothesis provides
$$Z_{(-\infty, \gamma_1)}\left(H_{2-\frac{1}{\varkappa}}[p']\right)=Z_{( \gamma_{n-2}, \infty)}
\left(H_{2-\frac{1}{\varkappa}}[p']\right)=1,$$
where $\gamma_1<\gamma_2<\ldots<\gamma_{n-2}$ are zeros of $p''.$ Since
 $$\xi_1<\gamma_1<\xi_2<\gamma_2<\ldots<\xi_{n-2}<\gamma_{n-2}<\xi_{n-1},$$ 
we have
\begin{equation}
\label{bbaa60}
Z_{(-\infty, \xi_1)}\left(H_{2-\frac{1}{\varkappa}}[p']\right)\leq 1, \ \ Z_{( \xi_{n-1}, \infty)}
\left(H_{2-\frac{1}{\varkappa}}[p']\right)\leq 1.
\end{equation}
If we apply Corollary \ref{Rolle} from Lemma~\ref{lemma8}, we obtain
\begin{equation}
\label{bbaa61}
Z_{(-\infty, \xi_1)}\left(H_{\varkappa}[p]\right) \leq 2, \ \ Z_{( \xi_{n-1}, \infty)}
\left(H_{\varkappa}[p]\right)\leq 2.
\end{equation}
According to Lemma~\ref{Interval.endpoints},  the number $Z_{I_k}\left(Q_{\varkappa}[p]\right) $ 
is odd if $I_k$ is an infinite interval of the first type. Thus, for $\frac{j-1}{j}\leq \varkappa<
\frac{j}{j+1}$ we obtain
\begin{equation}
\label{bbaa62}
Z_{I_1}\left(H_{\varkappa}[p]\right)=1, \ \ Z_{I_n}\left(H_{\varkappa}[p]\right)=1.
\end{equation}
Lemma~\ref{lemma15} is proved. $\square$
\begin{remark}\label{GrinfI}
It follows from (\ref{tsss90}), Theorem~\ref{Theorem4}, and Lemma~\ref{lemma15} that the 
function $M[p]$ decreases on the interval $(-\infty, \xi_1)$ from $\frac{n-1}{n}$ to $-\infty,$ 
and increases on the interval $(\xi_{n-1},\infty)$ from $-\infty$ to $\frac{n-1}{n}.$  
\end{remark}

\section{Accuracy of estimates in Theorem~\ref{Theorem3}} \label{Sec14}

Now we are going to discuss accuracy of estimates in Theorem~\ref{Theorem3}. 
The following fact could be derived from the first and the third statements of 
Theorem~\ref{Theorem4} immediately. This theorem provides accuracy of the estimation 
from below in Theorem~\ref{Theorem3}.

\begin{corol}\label{lowex}
Let $0<\varkappa<1/2.$ Assume that $p_{2n}(x)=x^{2n}+\ldots$ is a polynomial having no real 
roots, and all roots of its derivative $p_{2n}'$ are real and simple. Then 
$$Z_{\mathbb{R}}(\varkappa( p_{2n}')^2-p_{2n}p''_{2n})=2n-2.$$
\end{corol}

Denote by $T_n(x)=cos(n \arccos (x))$ the Chebyshev polynomial of the first kind of degree $n$. 
The fact that the bound from above cannot be improved is a corollary of the statement below.
\begin{theorem}\label{Theorem7}
 Let $0<\varkappa<1/2.$ For every $n=1,2,\ldots,$ and each $0<\varepsilon<\frac{1}{2}$ there exists 
 a constant $C=C(\varepsilon, n), \  -1<C<0,$ such that for every $0<\varkappa\leq \frac{1}{2}-\varepsilon$ 
 the following is true
\begin{equation}
\label{aaab1}
Z_{\mathbb{R}}\left(\varkappa-M[T_{2n}+C]\right)=Z_{\mathbb{R}}\left(\varkappa-
\frac{(T_{2n}(x)+C)T_{2n}''(x)}{(T_{2n}'(x))^2}\right)=4n-2.
\end{equation}
\end{theorem}

{\bf Proof of Theorem~\ref{Theorem7}.} First, consider the case of $n=1,$ that is of $T_2(x)+C=2x^2-1+C.$ 
If $|C|<1,$ both roots of the polynomial are real and simple. Since $(T_2(x)+C)' =4x,$ this polynomial 
generates two intervals $I_k$: $I_1=(-\infty, \ 0)$ and $I_2=(0, \ +\infty).$ Both intervals are intervals of 
the first type. Thus, the statement of Theorem~\ref{Theorem7} follows immediately from 
Remark~\ref{GrinfI}.

Now, for any $n=2,3, \ldots,$ we consider the polynomial
$$P_{2n}(x)=T_{2n}(x)-1=2^{2n-1}\prod_{k=0}^{2n-1}\left( x-cos
\left( \frac{k\pi}{2n} \right) \right)$$
\begin{equation}
\label{bbaa3}
=2^{2n-1}(x^2-1) \prod_{k=1}^{n-1}\left( x-\cos\left( \frac{k\pi}{n} \right) \right)^2.
\end{equation}

By (\ref{limit}) for any $k=1,\ldots, n-1$ we have
$$ \lim_{x\to \cos\left( \frac{k\pi}{n} \right) }M[P_{2n}](x)=\frac{1}{2}.$$ 
Let us fix any $0<\varepsilon<\frac{1}{2}$. There are numbers $\delta_k=\delta_k(\varepsilon)$ 
such that for all $x,\ |x-\cos\left( \frac{k\pi}{n}\right)|<\delta_k, \  k=1,2,\ldots, n-1,$ the following is true 
\begin{equation}
\label{bbaa4}
\left|M[P_{2n}](x)-\frac{1}{2}\right|<\frac{\varepsilon}{4}.
\end{equation}
Consider the family of polynomials $\{P_{2n}(x)+B, \ \ 0<B<1 \}$ on the compact set
\begin{equation*}
\label{bbaa15}
K=\cup_{k=1}^{n-1} \left(\left[\cos\left(\frac{k\pi}{n}\right)-\frac{\delta_k}{2}, 
\cos \left(\frac{k\pi}{n}\right)-\frac{\delta_k}{4}\right]\cup \left[\cos \left(\frac{k\pi}{n}\right)
+\frac{\delta_k}{4}, \cos \left(\frac{k\pi}{n}\right)+\frac{\delta_k}{2}\right]\right).
\end{equation*}
Obviously, on $K$ we have
\begin{equation}
\label{bbaa5}
M \left[P_{2n}+B\right](x)=\frac{(P_{2n}(x)+B)P_{2n}''(x)}{(P_{2n}(x))^2}
\rightrightarrows M[P_{2n}](x), \ \ B \to 0.
\end{equation}
Therefore, there exists $B_0=B_0(\varepsilon), \ 0<\ B_0<1/2,$ such that for all 
$0<B<B_0$ and for all $x \in K$ the inequality is true
\begin{equation}
\label{bbaa6}
\left|M \left[P_{2n}+B\right](x)-M[P_{2n}](x)\right|<\frac{\varepsilon}{4}.
\end{equation}
It follows from (\ref{bbaa4}) and (\ref{bbaa6}) that for all $0<B<B_0$ and for all $x \in K$ we have
\begin{equation}
\label{bbaa7}
\left|M \left[P_{2n}+B\right](x)-\frac{1}{2} \right|<\frac{\varepsilon}{2}.
\end{equation}
Thus, for all $0<B<B_0$ and for all $x \in K$ we obtain the estimation
\begin{equation}
\label{bbaa8}
M \left[P_{2n}+B\right](x)=\frac{(P_{2n}(x)+B)P_{2n}''(x)}{(P_{2n}(x))^2}>
\frac{1}{2}-\frac{\varepsilon}{2}>0.
\end{equation}
Note that 
$$(P_{2n}(x)+B)'=(T_{2n}(x)-1+B)'=T_{2n}'(x)=$$
$$2n\frac{\sin(2n\arccos x)}{\sin(\arccos x)}=2^{2n}n \prod_{k=1}^{2n-1}
\left( x-\cos\left( \frac{k\pi}{2n} \right) \right).$$
Since $0<B<\frac{1}{2},$ the roots of $P_{2n}(x)+B=T_{2n}(x)-1+B$ are all real and simple. 
Therefore, all the intervals $I_k, \ k=1,2,\ldots, 2n,$ are intervals of the first type, where
$$I_1=\left(-\infty, \ \cos\left( \frac{\pi}{2n}\right) \right), \ \  I_{2n}=\left( 
\cos\left( \frac{(2n-1)\pi}{2n}\right),\ \infty \right), $$
\begin{equation}
\label{bbaa12}
 \ I_{k}=\left(\cos\left( \frac{(k-1)\pi}{2n}\right), \ 
 \cos\left( \frac{k\pi}{2n}\right) \right), \ \ k=2,\ldots ,2n-1. 
\end{equation}
Making if necessary the values of $\delta_k, \  k=1,2, \ldots, n-1,$ smaller we can assume that
$$\left[\cos \left(\frac{k\pi}{n}\right)-\frac{\delta_k}{2}, \cos \left( \frac{k\pi}{n}\right)-
\frac{\delta_k}{4}\right] \subset I_{2k}, $$ 
$$\left[\cos \left(\frac{k\pi}{n}\right)+\frac{\delta_k}{4}, \cos \left(\frac{k\pi}{n}\right)+
\frac{\delta_k}{2}\right] \subset I_{2k+1}, \ \ k=1, 2,\ldots, n-1. $$
It follows from the second statements of Theorem~\ref{Theorem4}and (\ref{bbaa8}), that for each 
$0<\varkappa\leq \frac{1}{2}-\varepsilon$ on any finite interval $I_k, \  k=2,3,\ldots, 2n-1,$ the 
following is true
\begin{equation}
\label{bbaa9}
Z_{I_k}\left(\varkappa-\frac{(P_{2n}(x)+B)P_{2n}''(x)}{(P_{2n}(x))^2}\right)=2.
\end{equation}
Because all roots of $P_{2n}(x)+B$ are real and simple, by virtue of the fourth statement of 
Theorem~\ref{Theorem4} for each $0<\varkappa<\frac{1}{2}$ on each of the infinite intervals 
$I_k$ we have 
\begin{equation}
\label{bbaa10}
Z_{I_k}\left(\varkappa-\frac{(P_{2n}(x)+B)P_{2n}''(x)}{(P_{2n}(x))^2}\right)=1.
\end{equation}
Finally we obtain
\begin{equation}
\label{bbaa11}
Z_{\mathbb{R}}\left(\varkappa-\frac{(P_{2n}(x)+B)P_{2n}''(x)}{(P_{2n}(x))^2}\right)
=2(2n-2)+2=4n-2.
\end{equation}

Theorem~\ref{Theorem7} is proved.  $\square$ 

\section{Polynomials with only real zeros} \label{Sec15}

As we mentioned in the introductory part of this paper, in their paper \cite{Tyaglov.Atia.2021} 
M.~Tyaglov and M.J.~Atia found the upper and lower bounds for $Z_{\mathbb{R}}\left(H_{\varkappa}\right)[p]$ 
for all real $\varkappa$ when all roots of a polynomial $p$ are real. In Theorem \ref{TheoremTA} we described 
their result. In this section, we are going to refine their result in the following way. We discuss how the roots of the 
corresponding function $Q_{\varkappa}[p]$ can be distributed among the finite intervals $I_k.$ 

\begin{lemma}\label{lemma17} For all $a_1\geq 0, a_2\geq 0, \ldots, a_m\geq 0$ the following 
inequalities are valid
\begin{equation}
\label{bbaa23}
2\sum_{1\leq j<k \leq m} a_j a_k \ \leq \ (m-1) \sum_{k=1}^{m} a_k ^2, 
\end{equation}
\begin{equation}
\label{bbaa21}
\left(\sum_{k=1}^{m} a_k \right)^2 \ \leq \ m \sum_{k=1}^{m} a_k^2, 
\end{equation}
\begin{equation}
\label{bbaa22}
 (a_1+\ldots +a_s-a_{s+1}-\ldots -a_m)^2 \ \leq  \ \max(s, m-s)  
 \sum_{k=1}^{m} a_k^2 , \  \  1\leq s \leq m-1,
\end{equation}
and, if all $a_j$ are nonzero numbers for $j=1, 2, \ldots, m,$ then 
\begin{equation}
\label{bbaa222}
 (a_1+\ldots +a_s-a_{s+1}-\ldots -a_m)^2 \ <  \ \max(s, m-s)  
 \sum_{k=1}^{m} a_k^2 , \  \  1\leq s \leq m-1.
\end{equation}
\end{lemma}
{\bf Proof of Lemma~\ref{lemma17}}. The inequality (\ref{bbaa21}) is equivalent to the 
inequality (\ref{bbaa23}), because
$$\left(\sum_{k=1}^{m} a_k \right)^2=\sum_{k=1}^{m} a_k^2+2\sum_{1\leq j<k\leq m} a_j a_k. $$
Let us prove (\ref{bbaa23})
$$2 \sum_{1\leq j<k\leq m} a_j a_k  = \sum_{j=1}^{m}\sum_{k=j+1}^{m}(2a_j a_k)\leq
 \sum_{j=1}^{m}\sum_{k=j+1}^{m}(a_j^2+a_k^2)=\sum_{j=1}^{m}\left((m-j)a_j^2+
 \sum_{k=j+1}^{m}a_k^2\right)$$
$$=\sum_{j=1}^{m} (m-j)a_j^2+\sum_{k=2}^{m}\sum_{j=1}^{k-1}a_k^2=\sum_{j=1}^{m} 
(m-j)a_j^2+\sum_{k=2}^{m}(k-1)a_k^2$$
$$=(m-1)a_1^2+\sum_{k=2}^{m}\left((m-k)a_k^2+(k-1)a_k^2\right)=(m-1) \sum_{k=1}^{m} a_k^2.$$
The inequalities (\ref{bbaa23}) and (\ref{bbaa21}) are proved.

Let us fix $1\leq s \leq m-1$ and prove (\ref{bbaa22}). We have
$$\left(\sum_{k=1}^{s}a_k -\sum_{k=s+1}^{m}a_k \right)^2 =
\left(\sum_{k=1}^{s}a_k \right)^2+\left(\sum_{k=s+1}^{m}a_k \right)^2-
2\left(\sum_{k=1}^{s}a_k \right) \cdot \left(\sum_{k=s+1}^{m}a_k \right)$$
$$\leq \left(\sum_{k=1}^{s}a_k \right)^2+\left(\sum_{k=s+1}^{m}a_k \right)^2,$$
where the inequality sign is strict if all numbers are nonzero.
Applying (\ref{bbaa21}) to each of the terms in the right-hand side of the above inequality we obtain
$$\left(\sum_{k=1}^{s}a_k -\sum_{k=s+1}^{m}a_k \right)^2 \ \leq \ s \sum_{k=1}^{s} 
a_k^2 +(m-s) \sum_{k=s+1}^{m} a_k^2 \leq \max(s, m-s)  \sum_{k=1}^{m} a_k^2.$$
Lemma~\ref{lemma17} is proved. $\square$

The following statement about maximum values of the function $M[p]$ on a finite interval $I_k$ in the case 
when all intervals $I_k$ are intervals of the first type, will be obtained as a corollary of the above lemma.
\begin{theorem}\label{Theorem8} Assume that a polynomial $p$  has only 
real and simple zeros. Then for all $x\in I_s \cup I_{n-s+1}=(\xi_{s-1}, \ \xi_s)\cup (\xi_{n-s}, \ \xi_{n-s+1}), 
\ \ s=2,3,\ldots, \lfloor\frac{n+1}{2}\rfloor,  $ the following estimation is valid
\begin{equation}
\label{bbaa24}
M[p](x)=\frac{p(x) p''(x)}{(p'(x))^2} <\frac{n-s}{n-s+1}.  
\end{equation}
\end{theorem}
{\bf Proof of Theorem~\ref{Theorem8}}. Suppose that $p(x)=(x-x_1)(x-x_2)\cdots (x-x_n).$  We will assume 
that $x\neq x_k, \ \ k=1,2,\ldots, n.$  Since all  roots of the polynomial $p$ are simple, if $x=x_k$, then 
$M_p(x_k)=0$. In this case the estimation (\ref{bbaa24}) is obvious.

Let $x\neq x_k.$ We have
$$\frac{p(x)}{p'(x)} =\left(\frac{p'(x)}{p(x)}\right)^{-1}=\left(\sum_{k=1}^n \frac{1}{x-x_k}\right)^{-1}.$$
Therefore,
\begin{equation}
\label{bbaa25}
\left(\frac{p(x)}{p'(x)}\right)'=1-\frac{p(x) p''(x)}{(p'(x))^2} =\frac{\sum_{k=1}^n 
\left(\frac{1}{x-x_k}\right)^2}{\left(\sum_{k=1}^n \frac{1}{x-x_k}\right)^2}.
\end{equation}
As before we will denote by $\xi$ with indices the roots of the derivative $p'.$
Since all roots of $p$ and $p'$ are real and simple, we have
\begin{equation}\label{ABC1}
x_1<\xi_1<x_2<\xi_2<\ldots<\xi_{n-1}<x_n.
\end{equation}
So, if $x \in I_s,  \ s=2,3,\ldots, \lfloor\frac{n+1}{2}\rfloor,$ then $$x_{s-1}<\xi_{s-1}<x<\xi_s 
<x_{s+1}.$$ Therefore, 
$$ x-x_{k}>0, \ \ \mbox{when} \ \ k=1,2,\ldots, s-1, \ \ \mbox{and}$$  
\begin{equation}\label{ABC2}
 x-x_{k}<0, \ \ \mbox{when} \ \ k=s+1,s+2,\ldots, n,
\end{equation}
for all $x \in I_s.$
Similar if  $x \in I_{n-s+1}, \ \ s=2,3,\ldots, \lfloor\frac{n+1}{2}\rfloor,$ then $$x_{n-s}<\xi_{n-s}<
x<\xi_{n-s+1}<x_{n-s+2}.$$ Thus,
$$ x-x_{k}>0, \ \ \mbox{when} \ \ k=1,2,\ldots, n-s, \ \ \mbox{and}$$  
\begin{equation}\label{ABC3}
 x-x_{k}<0, \ \ \mbox{when} \ \ k=n-s+2,\ldots, n,
\end{equation}
for all $x \in I_{n-s+1}.$
Note that for $s< \lfloor\frac{n+1}{2}\rfloor$ we have $s\leq n-s+1.$ Applying (\ref{bbaa222}) to (\ref{bbaa25}), we obtain
\begin{equation}
\label{bbaa26}
1-\frac{p(x) p''(x)}{(p'(x))^2} =\frac{\sum_{k=1}^n \left(\frac{1}{x-x_k}\right)^2}{\left(\sum_{k=1}^n \frac{1}{x-x_k}\right)^2}
>\frac{1}{n-s+1}.
\end{equation}
The statement (\ref{bbaa24}) of Theorem~\ref{Theorem8} follows from the last inequality.
 Theorem~\ref{Theorem8} is proved.  $\square$

\begin{corol}\label{C1T8} Assume that a polynomial $p$  has only real 
and simple zeros. If  $$\frac{n-2}{n-1}\leq \varkappa<\frac{n-1}{n},$$ then for all finite intervals $I_k$ 
the following is true
\begin{equation}
\label{bbaa27}
Z_{I_k}(H_{\varkappa}[p])=Z_{I_k}\left(Q_{\varkappa}[p]\right) =0.
\end{equation}
\end{corol}
The more general statement below contains Corollary~\ref{C1T8} as a partial case.
\begin{corol}\label{C2T8}
Assume that a polynomial $p$  has only real and simple zeros. 
Let for some $j= \lfloor\frac{n}{2}\rfloor+1,\ldots, n-1,$ the inequality 
$$\frac{j-1}{j}\leq \varkappa<\frac{j}{j+1}$$ 
holds.  Then for any interval $I_k$, $k=n-j+1,\ldots,j,$ the following is true
\begin{equation}
\label{bbaa27}
Z_{I_k}(H_{\varkappa}[p])=Z_{I_k}\left(Q_{\varkappa}[p]\right) =0.
\end{equation}
\end{corol}

The theorem below follows from Theorem~\ref{Theorem5} and Corollary~\ref{C2T8}.
\begin{theorem}\label{Theorem9} Assume that a polynomial $p$  
has only real and simple zeros. Let for some $j=1,2,\ldots, n-1$ the inequality $$\frac{j-1}{j}\leq 
\varkappa<\frac{j}{j+1}$$ holds. Then 
\begin{enumerate}
\item { \it If $j=\lfloor\frac{n}{2}\rfloor+1,\ldots, n-1,$ then} 
\begin{itemize}
\item {\it For any finite interval $I_k$, where $k=n-j+1,\ldots,j,$ the estimation is valid}
$$Z_{I_k}(H_{\varkappa}[p])=Z_{I_k}\left(Q_{\varkappa}[p]\right) =0.$$
\item {\it  For any finite interval $I_k$, where $k=2,\ldots, n-j,$ or $k=j+1,\ldots, n-1,$ the estimation is valid }
$$Z_{I_k}(H_{\varkappa}[p])=Z_{I_k}\left(Q_{\varkappa}[p]\right) =0 \ \mbox{or} \  2.$$
\end{itemize}
\item { \it If $j=1, 2,\ldots, \lfloor\frac{n-1}{2}\rfloor,$ then for any finite interval $I_k$ the estimation is valid}
$$Z_{I_k}(H_{\varkappa}[p])=Z_{I_k}\left(Q_{\varkappa}[p]\right) =0 \ \mbox{or} \ 2.$$
\end{enumerate}
\end{theorem}
The following result shows  that Theorem~\ref{Theorem9} cannot be improved.
\begin{theorem}\label{Theorem10}
For each $\varepsilon>0$ there exists a polynomial 
\begin{equation}
\label{abcd1}
p_{\varepsilon}(x)=(x-x_1(\varepsilon))(x-x_2(\varepsilon))\cdot \ldots \cdot(x-x_n(\varepsilon)),
\end{equation}
and a system of $n-2$ points $y_k (\varepsilon) \in I_k, \ \ k=2,3,\ldots, n-1, $ such that
\begin{equation}
\label{abcd2}
p_{\varepsilon}(y_k (\varepsilon))>\frac{k-1}{k}-\varepsilon.
\end{equation}
\end{theorem}

{\bf Proof of Theorem~\ref{Theorem10}.} We will construct the polynomial $p_{\varepsilon}$ using 
inductive process. We start with the polynomial
\begin{equation}
\label{abcd4}
p_0(x)=x^{n-1}(x-1).
\end{equation}
Its derivative has the form
\begin{equation}
\label{bacd1}
p'_0(x)=nx^{n-2}\left(x-\frac{n-1}{n}\right).
\end{equation}
By virtue of (\ref{limit}) and definition (\ref{ttss27}) of the function $M[p]$ we have
\begin{equation}
\label{abcd5}
\lim_{x\to 0} M[p_0](x)=\frac{n-2}{n-1}.
\end{equation}
Let us fix any (small enough) $\varepsilon >0.$  It follows from (\ref{abcd5}) that 
there is $0<\delta_1<1$ such that
\begin{equation}
\label{abcd6}
 M[p_0](x)>\frac{n-2}{n-1}-\frac{\varepsilon}{n}, \ \ \mbox{whenever} \ |x|<\delta_1.
\end{equation}
Denote by
\begin{equation}
\label{abcd7}
p_1(x)=x^{n-2}(x-x_{n-1})(x-1),
\end{equation}
where $0<x_{n-1}<1$ will be chosen later.
The derivative of the above polynomial will have the form
\begin{equation}
\label{bacd2}
p_1'(x)=x^{n-3}(x-\xi_{n-2}^1)(x-\xi_{n-1}^1),
\end{equation}
where $0<\xi_{n-2}^1<x_1<\xi_{n-1}^1<1$.
Note that on the compact set $\frac{\delta_1}{4} \leq x \leq \frac{\delta_1}{2}$ the following is true
\begin{equation}
\label{abcd8}
M[p_1](x) \rightrightarrows M[ p_0](x), \ \  \mbox{when}  \  \ x_{n-1}\rightarrow 0,
\end{equation}
and
\begin{equation}
\label{bacd3}
\xi_{n-2}^1\rightarrow 0 \ \ \mbox{and}  \  \ \xi_{n-1}^1 \rightarrow \frac{n-1}{n}, \ \ 
\mbox{when}  \  \ x_{n-1}\rightarrow 0.
\end{equation}
Thus, we can find a number $x_{n-1}=x_{n-1}(\varepsilon)$ such that 
\begin{equation}
\label{bacd4}
\xi_{n-2}^1< \frac{\delta_1}{4} < \frac{\delta_1}{2}<\xi_{n-1}^1, 
\end{equation}
and for all $x\in \left[\frac{\delta_1}{4}, \frac{\delta_1}{2} \right]$ we have
\begin{equation}
\label{abcd9}
M[p_1](x)>M[p_0](x)-\frac{\varepsilon}{n},
\end{equation}
where
\begin{equation}
\label{abcd10}
p_1(x)=x^{n-2}(x-x_{n-1}(\varepsilon))(x-1).
\end{equation}

Let us fix in (\ref{abcd10}) $x_{n-1}=x_{n-1}(\varepsilon)$ that provides (\ref{bacd4}) and 
(\ref{abcd9}). It follows from (\ref{abcd6}) and (\ref{abcd9}) that
\begin{equation}
\label{abcd11}
M[p_1](x)>\frac{n-2}{n-1}-\frac{2\varepsilon}{n},
\end{equation}
for all $x\in \left[\frac{\delta_1}{4}, \frac{\delta_1}{2} \right].$

For the next step we consider (\ref{abcd10}) as the initial polynomial. Its derivative has the 
form (\ref{bacd2}). By virtue of (\ref{limit})
\begin{equation}
\label{abcd12}
\lim_{x\to 0} M[ p_1](x)=\frac{n-3}{n-2}.
\end{equation}
By (\ref{abcd12}) there is $0<\delta_2<\frac{\delta_1}{2}$ such that
\begin{equation}
\label{abcd13}
 M[p_1](x)>\frac{n-3}{n-2}-\frac{\varepsilon}{n}, \ \ \mbox{whenever} \ |x|<\delta_2.
\end{equation}
Denote by
\begin{equation}
\label{abcd14}
p_2(x)=x^{n-3}(x-x_{n-2})(x-x_{n-1}(\varepsilon))(x-1),
\end{equation}
where $0<x_{n-2}<x_{n-1}(\varepsilon)$ will be chosen later.
The derivative of the above polynomial will have the form
\begin{equation}
\label{bacd5}
p_1'(x)=x^{n-4}(x-\xi_{n-3}^2)(x-\xi_{n-2}^2)(x-\xi_{n-1}^2),
\end{equation}
where $0<\xi_{n-3}^2<x_{n-2}<\xi_{n-2}^2<x_{n-1}(\varepsilon)<\xi_{n-1}^2<1$.
Note that on the compact set $\left[\frac{\delta_2}{4} , \frac{\delta_2}{2}\right] 
\cup \left[\frac{\delta_1}{4} , \frac{\delta_1}{2}\right]$ the following is true
\begin{equation}
\label{abcd15}
M[p_2](x) \rightrightarrows M[p_1](x), \ \   \mbox{when}  \  \  x_{n-2}\rightarrow 0,
\end{equation}
and
\begin{equation}
\label{bacd6}
\xi_{n-3}^2\rightarrow 0, \  \ \xi_{n-2}^2 \rightarrow \xi_{n-2}^1,\  \  \xi_{n-1}^2 \rightarrow \xi_{n-1}^1,
\ \ \mbox{when}  \  \ x_{n-2}\rightarrow 0.
\end{equation}
Taking into account (\ref{bacd4}) and the fact that $0<\delta_2<\frac{\delta_1}{2},$ we can find a number 
$x_{n-2}=x_{n-2}(\varepsilon)$ such that 
\begin{equation}
\label{bacd7}
\xi_{n-3}^2< \frac{\delta_2}{4} < \frac{\delta_2}{2}<\xi_{n-2}^2< \frac{\delta_1}{4} 
< \frac{\delta_1}{2}<\xi_{n-1}^2, 
\end{equation}
and such that for all $x\in \left[\frac{\delta_2}{4} , \frac{\delta_2}{2}\right] \cup \left[\frac{\delta_1}{4} , 
\frac{\delta_1}{2}\right]$ we have
\begin{equation}
\label{abcd16}
M[p_2](x)>M[p_1](x)-\frac{\varepsilon}{n},
\end{equation}
where
\begin{equation}
\label{abcd17}
p_2(x)=x^{n-3}(x-x_{n-2}(\varepsilon))(x-x_{n-1}(\varepsilon))(x-1).
\end{equation}
It follows from (\ref{abcd13}) and (\ref{abcd16}) that
\begin{equation}
\label{abcd18}
M[p_2](x)>\frac{n-3}{n-2}-\frac{2\varepsilon}{n},
\end{equation}
for all $x\in \left[\frac{\delta_2}{4}, \frac{\delta_2}{2} \right],$
and from (\ref{abcd11}) and (\ref{abcd16}) that
\begin{equation}
\label{abcd19}
M[p_2](x)>\frac{n-2}{n-1}-\frac{3\varepsilon}{n},
\end{equation}
for all $x\in \left[\frac{\delta_1}{4}, \frac{\delta_1}{2} \right].$

If we continue this process, finally we will obtain a polynomial
$$p_{\varepsilon}(x)=(x-x_1(\varepsilon))(x-x_2(\varepsilon)) \cdot \ldots \cdot (x-x_n(\varepsilon)),$$
and the system of intervals 
$$\bigcup_{k=1}^{n-2} \left[\frac{\delta_k}{4}, \frac{\delta_k}{2} \right],$$
where
\begin{equation}
\label{abcd20} 
\left[\frac{\delta_k}{4}, \frac{\delta_k}{2} \right] \subset  I_{n-k}, \ \ k=1, 2, \ldots, n-2,
\end{equation}
such that for each $k=1,2,\ldots, n-2,$ the inequalities are true
\begin{equation}
\label{abcd21} 
M[p_{\varepsilon}](x)>\frac{n-k-1}{n-k}-\frac{(n-k+1)\varepsilon}{n}, \ \ x \in 
\left[\frac{\delta_k}{4}, \frac{\delta_k}{2} \right]. 
\end{equation}
The statement of Theorem~\ref{Theorem10} follows from (\ref{abcd20}) and (\ref{abcd21}). 
Theorem~\ref{Theorem10} is proved.

\section{Intervals $I_k$ of the second type} \label{Sec16}

Return now to the general case. The statement below gives an estimation for the number 
of real zeros of a polynomial $H_{\varkappa}[p]$ on infinite intervals $I_k$ of the second type.

\begin{lemma}\label{lemma19} Assume that $I_k$ is an infinite interval of the second type 
for a polynomial $p$ whose derivative $p'$ has only real and simple zeros. Then there exists a number 
$$C=\min_{x \in I_k} \frac{p(x)p''(x)}{(p'(x))^2}, \ \  \frac{1}{2}<C<\frac{n-1}{n},$$
such that
\begin{itemize}
\item {\it For every $C\leq \varkappa <\frac{n-1}{n}$ the following estimation is valid}
\begin{equation}
\label{abd1}
Z_{I_k}(H_{\varkappa}[p])=Z_{I_k}(Q_{\varkappa}[p])=2.
\end{equation}
\item {\it For every $\frac{1}{2} \leq \varkappa <C$ the following estimation is valid}
\begin{equation}
\label{abd2}
Z_{I_k}(H_{\varkappa}[p])=Z_{I_k}(Q_{\varkappa}[p])=0.
\end{equation}
\end{itemize}
\end{lemma}
{\bf Proof of Lemma~\ref{lemma19}.} If $\frac{1}{2}\leq \varkappa<\frac{n-1}{n},$ then 
$0\leq 2-\frac{1}{\varkappa}<\frac{n-2}{n-1}.$ Since all roots of $p'$ are real and simple, the 
intervals $(-\infty, \gamma_1),$ $(\gamma_{n-2}, +\infty),$ where $\gamma_1$ and 
$\gamma_{n-2}$ are the smallest and the biggest roots of the polynomial $p'',$ are intervals 
of the first type for the polynomial $p',$ and by Theorem~\ref{Theorem6}.
$$Z_{(-\infty, \gamma_1)}\left(H_{2-\frac{1}{\varkappa}}[p']\right)=
Z_{( \gamma_{n-2}, +\infty)}\left(H_{2-\frac{1}{\varkappa}}[p']\right)=1.$$
Since $\xi_1<\gamma_1<\ldots<\gamma_{n-2}<\xi_{n-1},$ we have
$$Z_{(-\infty, \xi_1)}\left(H_{2-\frac{1}{\varkappa}}[p']\right)\leq 1, \ \ 
Z_{( \xi_{n-1}, +\infty)}\left(H_{2-\frac{1}{\varkappa}}[p']\right)\leq 1.$$
By virtue of Corollary~\ref{Rolle} from Lemma~\ref{lemma8}, if $I_k$ is an infinite interval of 
the second type, then
$$Z_{I_k}\left(H_{\varkappa}[p]\right)\leq 2.$$
Taking into account the fourth statement of Lemma~\ref{Infint}, finally we obtain
\begin{equation}
\label{abd3}
Z_{I_k}\left(H_{\varkappa}(p)\right) \ =\ 0 \ \ \mbox{or} \  2.
\end{equation}
Since all roots of the polynomial $p'$ are real, the polynomial $p$ satisfies the assumptions 
of Theorem~\ref{conjectureSh2}. In the same way as in Theorem~\ref{conjectureSh2} we can 
show that (\ref{abd10}) is true, that is there exists such a value $x_0$ that for $|x|>x_0>0$ 
$$M[p](x)-\frac{n-1}{n}<0.$$ 
So, the graph of the function $y=M[p](x)$ approaches its asymptote $y=\frac{n-1}{n}$ from below. 

 If $I_k$ is an infinite interval of the second type, then its finite endpoint is a wrong point 
 (Remark~\ref{Intend}). We will denote it just by $\xi$ without indices. By (\ref{IRW}) and (\ref{ILW}),
\begin{equation}
\label{abd5}
\lim_{x\to \xi} M[p](x) \ = +\infty.
\end{equation}
Denote by
\begin{equation}
\label{abd13}
C=\min_{x \in I_1\cup I_{n}}M[p](x).
 \end{equation}
It follows from (\ref{abd5}) and (\ref{abd10}) that
\begin{equation}
\label{abd11}
C<\frac{n-1}{n},
\end{equation}
and from the third statement of Theorem~\ref{Theorem4} that
\begin{equation}
\label{abd12}
C > \frac{1}{2}.
\end{equation}
Lemma~\ref{lemma19} is proved. $\square$

Let $I_k$ be a finite interval of the second type. As we noticed in Remark~\ref{Intend}, one of the 
endpoints of $I_k$ is right while the other is wrong. It follows from (\ref{IRR}) -- (\ref{ILW}) 
that for any real $\varkappa$ the equation (\ref{mainequation}) 
$$M[p](x) =\varkappa$$
has solutions in each finite interval $I_k$ of the second type. Putting together (\ref{II2}), (\ref{sstt13}), 
and the first statement of Theorem~\ref{Theorem4}, we conclude that the equation (\ref{mainequation}) 
has the unique solution for all $\varkappa \in (-\infty, \frac{1}{2})\cup (\frac{n-1}{n}, \infty).$ It means 
that on corresponding intervals inside $I_k$ the function $M[p]$ is monotone. 

There may be an assumption that the function $M[p]$ is monotone on any finite interval $I_k$ of the 
second type. If we use Corollary~\ref{Rolle} from Lemma~\ref{lemma8}, in the same way as above, we 
can easily obtain that for $\frac{1}{2}\leq \varkappa < \frac{n-1}{n}$ on each interval $I_k$ of the second 
type
\begin{equation}\label{abd50}
Z_{I_k}\left(Q_{\varkappa}[p] \right)=Z_{I_k}\left(H_{\varkappa}[p] \right)\leq3.
\end{equation}

The following example shows that the estimation (\ref{abd50}) cannot be improved, that is there is a 
polynomial $p$ all whose real roots are simple, and the roots of its derivative $p'$ are real and simple, 
such that $Z_{I_k}\left(Q_{\varkappa}[p] \right)=3$ on some of the finite intervals $I_k$ of the second 
type. 

\begin{example} Consider a polynomial $p(x)=x^2(x-1)(x-2)(x+10)+0.1.$ Its derivative
$p'(x)=5x^4+28x^3-84x^2+40x$ has four real roots
$$\xi_1\approx -7.865, \ \xi_2=0, \  \xi_3 \approx 0.617 , \ \xi_4 \approx1.648.  $$
The equation
$$M[p](x)=\frac{2}{3}$$
has three solutions in the interval $(\xi_1,0).$
\end{example}
Indeed, it is easy to check that 
$$\lim_{x\to 0^{-}}M[p](x)=+\infty, \  
\lim_{x\to \xi_1^{+}}M[p](x)=-\infty, $$ 
$$ M[p](-0.275)\approx 0.642, \   
M[p](-1.35)\approx 0.684.$$

\setcounter{equation}{0}



\end{document}